\documentclass[aop,secthm,seceqn,MSNbibl,dvips]{arximspdf}

\doi{10.1214/11-AOP661}
\volume{40}
\issue{4}
\pubyear{2012}
\firstpage{1759}
\lastpage{1794}

\makeatletter
\newtheorem{hypo}[thm]{Hypothesis}
\newtheorem{lemma}[thm]{Lemma}
\newproclaim{remark}[thm]{Remark}
\newproclaim{defin}[thm]{Definition}
\renewcommand{\l}{l}
\makeatother

\begin{document}
\begin{frontmatter}

\title{The stochastic reflection problem on an infinite dimensional
convex set and BV functions in~a~Gelfand triple\thanksref{TITL1}}
\runtitle{Reflected OU-processes, BV functions}
\thankstext{TITL1}{Supported by 973 project, NSFC, key Lab of CAS, the DFG through IRTG 1132 and CRC 701 and the I. Newton Institute, Cambridge, UK.}

\begin{aug}
\author[A]{\fnms{Michael} \snm{R\"{o}ckner}\corref{}\ead[label=e1]{roeckner@mathematik.uni-bielefeld.de}},
\author[B]{\fnms{Rong-Chan} \snm{Zhu}\ead[label=e2]{zhurongchan@126.com}}
and
\author[C]{\fnms{Xiang-Chan} \snm{Zhu}\ead[label=e3]{zhuxiangchan@126.com}}
\runauthor{M. R\"{o}ckner, R.-C. Zhu and X.-C. Zhu}
\affiliation{University of Bielefeld, Chinese Academy of Sciences and
University of Bielefeld, and Peking University and University of Bielefeld}
\address[A]{M. R\"{o}ckner\\
Department of Mathematics\\
University of Bielefeld\\
D-33615 Bielefeld\\
Germany\\
\printead{e1}} 
\address[B]{R.-C. Zhu\\
Institute of Applied Mathematics\\
Academy of Mathematics and Systems\\
\quad  Science\\
Chinese Academy of Sciences\\
Beijing 100190\\
China\\
and\\
Department of Mathematics\\
University of Bielefeld\\
D-33615 Bielefeld\\
Germany\\
\printead{e2}}
\address[C]{X.-C. Zhu\\
School of Mathematical Sciences\\
Peking University\\
Beijing 100871\\
China\\
and\\
Department of Mathematics\\
University of Bielefeld\\
D-33615 Bielefeld\\
Germany\\
\printead{e3}}
\end{aug}

\received{\smonth{9} \syear{2010}}
\revised{\smonth{2} \syear{2011}}

%
\begin{abstract}
In this paper, we introduce a~definition of BV functions in a~Gelfand
triple which is an extension of the definition of BV functions in
[\textit{Atti Accad. Naz. Lincei Cl. Sci. Fis. Mat. Natur. Rend.
Lincei (9)
Mat. Appl.} \textbf{21} (2010) 405--414]
by using Dirichlet form theory. By this definition, we can consider the
stochastic reflection problem associated with a~self-adjoint operator
$A$ and a~cylindrical Wiener process on a~convex set $\Gamma$ in a~Hilbert space~$H$. We prove the existence and uniqueness of a~strong
solution of this problem when $\Gamma$ is a~regular convex set. The
result is also extended to the nonsymmetric case. Finally, we extend
our results to the case when $\Gamma=K_\alpha$, where $K_\alpha=\{
f\in
L^2 (0,1)|f\geq-\alpha\},\alpha\geq0$.
\end{abstract}

%
\begin{keyword}[class=AMS]
\kwd{60GXX}
\kwd{31C25}
\kwd{60G60}
\kwd{26A45}.
\end{keyword}
\begin{keyword}
\kwd{Dirichlet forms}
\kwd{stochastic reflection problems}
\kwd{BV function}
\kwd{Gelfand triples}
\kwd{integration by parts formula in infinite dimensions}.
\end{keyword}

\end{frontmatter}

\section{Introduction}\label{sec1}
A definition of BV functions in abstract Wiener spaces has been given
by Fukushima in~\cite{12}, Fukushima and Hino in~\cite{13},
based upon Dirichlet form theory. In this paper, we introduce BV
functions in a~Gelfand triple, which is an extension of BV functions in
a Hilbert space defined in~\cite{1}.
Here, we use a~version of the Riesz--Markov representation theorem in
infinite dimensions proved by Fukushima using the quasi-regularity
of the Dirichlet form (see~\cite{17}) to give a~characterization of BV
functions.

In this paper, we consider the Dirichlet form
\[
\mathcal{E}^\rho(u,v)=\frac{1}{2}\int_H\langle Du,Dv\rangle\rho
(z)\mu(dz)
\]
(where~$\mu$ is a~Gaussian measure in~$H$ and $\rho$ is a~BV function)
and its associated process. By using BV functions, we obtain a~Skorohod-type representation for the associated process, if $\rho
=I_\Gamma$ and $\Gamma$ is a~convex set.

As a~consequence of these results, we can solve the following
stochastic differential inclusion in the Hilbert space~$H$:
%
\begin{equation}\label{eq1.1}
\cases{
dX(t)+\bigl(AX(t)+N_\Gamma(X(t))\bigr)\,dt\ni \,dW(t),\cr
X(0)=x,
}
\end{equation}
where our solution is strong (in the probabilistic sense), if $\Gamma$
is regular.
Here $A\dvtx D(A)\subset H\rightarrow H$ is a~self-adjoint operator.
$N_\Gamma(x)$ is the normal cone to $\Gamma$ at $x$ and $W(t)$ is a~cylindrical Wiener process in~$H$. The precise meaning of the above
inclusion will be defined in Section~\ref{sec5.2}. The solution to~(\ref{eq1.1}) is
called distorted (if $\rho=I_\Gamma$, reflected) Ornslein--Uhlenbek
(OU for short)-process.

Equation~(\ref{eq1.1}) was first studied (strongly solved) in~\cite{19}, when
$H=L^2(0,1)$, $A$ is the Laplace operator with Dirichlet or Neumann
boundary conditions and $\Gamma$ is the convex set of all nonnegative
functions of $L^2(0,1)$; see also~\cite{28}. In~\cite{6}, the authors
study the situation when $\Gamma$ is a~regular convex set with
nonempty interior. They get precise information about the corresponding
Kolmogorov operator, but did not construct a~strong solution to~(\ref{eq1.1}).

In this paper, we consider a~convex set $\Gamma$. If $\Gamma$ is a~regular convex set, we show that $I_\Gamma$ is a~BV-function and thus
obtain existence and uniqueness results for~(\ref{eq1.1}). By a~modification of
\cite{12} and using~\cite{7}, we obtain the existence of an (in the
probabilistic sense) weak solution to~(\ref{eq1.1}). Then, we prove pathwise
uniqueness. Thus, by a~version of the Yamada--Watanabe Theorem (see
\cite{15}), we deduce that~(\ref{eq1.1}) has a~unique strong solution. We also
consider the case when $\Gamma=K_\alpha$, where $K_\alpha=\{f\in L^2
(0,1)|f\geq-\alpha\},\alpha\geq0,$ and prove our result about
Skorohod-type representation and that $I_{K_\alpha}$ is a~BV function,
if $\alpha>0$.

The solution of the reflection problem is based on an integration by
parts formula. The connection to BV functions is given in Theorem~\ref{thm3.1}
below, which is a~key result of this paper. It asserts that the
integration by parts formula for $\rho\cdot\mu$ gives a~characterization of BV functions
$\rho$, in the case where~$\mu$ is a~Gaussian measure. This is an
extension of the characterization of BV functions in finite dimension.
But an integration by parts formula is in fact enough for the
reflection problem. This we show in Section~\ref{sec6}, exploiting the beautiful
integration by parts
formula for $K_\alpha,\alpha\geq0$, proved in~\cite{28}, which in
case $\alpha=0$, that is, $K_0=\{f\in L^2(0,1)\dvtx f\geq0\}$, is with respect
to a~non-Gaussian measure, namely a~Bessel bridge.
Theorem~\ref{thm3.1} applies to prove that $I_{K_\alpha}$ is a~BV function, but
only if $\alpha>0$.

This paper is organized as follows. In Section~\ref{sec2}, we consider the
Dirichlet form
and its associated distorted OU-process. We introduce BV functions in
Section~\ref{sec3}, by which we can get the Skorohod type representation for the
OU-process. In Section~\ref{sec4}, we analyze the reflected OU-process. In
Section~\ref{sec5}, we get the existence and uniqueness of the solution for
(\ref{eq1.1}) if $\Gamma$ is a~regular convex set. We also extend these
results to the nonsymmetric case. In Section~\ref{sec6}, we consider the case
when $\Gamma=K_\alpha$, where $K_\alpha=\{f\in L^2 (0,1)|f\geq
-\alpha\},\allowbreak\alpha\geq0$.

Finally, we would like to mention that apart from contributing to
develop the theory of BV functions on infinite dimensional spaces, one
main motivation of this paper is to provide the probabilistic
counterpart to
results in~\cite{6} and~\cite{7}, by exploiting Dirichlet form theory
and its associated potential theory.

\section{The Dirichlet form and the associated distorted OU-process}\label{sec2}
Let~$H$ be a~real separable Hilbert space (with scalar product $\langle
\cdot,\cdot\rangle$ and norm denoted by $|\cdot|$). We denote its
Borel $\sigma$-algebra by $\mathcal{B}(H)$. Assume that:

\begin{hypo}\label{hypo2.1}
$A\dvtx D(A)\subset H\rightarrow H$ is a~linear
self-adjoint operator on~$H$ such that $\langle Ax,x\rangle\geq\delta
|x|^2\  \forall x\in D(A)$ for some $\delta>0$ and $A^{-1}$
is of trace class.
\end{hypo}

Since $A^{-1}$ is trace class, there exists
an orthonormal basis $\{e_j\}$ in~$H$ consisting of eigen-functions for
$A$ with corresponding eigenvalues $\alpha_j\in\mathbb{R}$, $j\in
\mathbb{N}$, that is,
\[
Ae_j=\alpha_je_j,\qquad j\in\mathbb{N}.
\]
Then $\alpha_j\geq\delta$ for all $j\in\mathbb{N}$.

Below $D\varphi\dvtx H\rightarrow H$ denotes the Fr\'{e}chet-derivative of
a function $\varphi\dvtx H\rightarrow\mathbb{R}$. By $C_b^1(H)$ we shall
denote the set of all bounded differentiable functions with continuous
and bounded derivatives. For $K\subset H$, the space $C^1_b(K)$ is
defined as the space of restrictions of all functions in $C_b^1(H)$ to
the subset $K$.~$\mu$ will denote the Gaussian measure in~$H$ with
mean $0$ and covariance operator
\[
Q:=\tfrac{1}{2}A^{-1}.
\]
Since $A$ is strictly positive,~$\mu$ is nondegenerate and has full
topological support. Let $L^p(H,\mu), p\in[1,\infty]$, denote the
corresponding real $L^p$-spaces equipped with the usual norms $\|\cdot
\|_p$. We set
\[
\lambda_j:=\frac{1}{2\alpha_j}\qquad  \forall j\in\mathbb{N}
\]
so that
\[
Qe_j=\lambda_je_j \qquad \forall j\in\mathbb{N}.
\]
For $\rho\in L^1_+(H,\mu)$ we consider
\[
\mathcal{E}^\rho(u,v)=\frac{1}{2}\int_H\langle Du,Dv\rangle\rho
(z)\mu(dz),\qquad u,v\in C_b^1(F),\vadjust{\goodbreak}
\]
where $F:=\operatorname{Supp}[\rho\cdot\mu]$ and $ L^1_+(H,\mu)$ denotes the set
of all nonnegative elements in $ L^1(H,\mu)$. Let $\operatorname{QR}(H)$ be the set
of all functions $\rho\in L^1_+(H,\mu)$ such that $(\mathcal{E}^\rho
, C_b^1(F))$ is closable on $L^2(F,\rho\cdotp\mu$). Its closure is
denoted by $(\mathcal{E}^\rho, \mathcal{F}^\rho)$. We denote by
$\mathcal{F}_e^\rho$ the extended Dirichlet space of $(\mathcal
{E}^\rho, \mathcal{F}^\rho)$, that is, $u\in\mathcal{F}_e^\rho$
if and only if $|u|<\infty$ $\rho\cdot\mu$-a.e. and there
exists a~sequence $\{u_n\}$ in~$\mathcal{F^\rho}$ such that $\mathcal
{E}^\rho(u_m-u_n,u_m-u_n)\rightarrow0$ as $n\geq m\rightarrow\infty
$ and $u_n\rightarrow u $ $  \rho\cdot\mu$-a.e. as
$n\rightarrow\infty$.

\begin{thm}\label{thm2.2}
Let $\rho\in \operatorname{QR}(H)$. Then $(\mathcal{E}^\rho,
\mathcal{F}^\rho)$ is a~quasi-regular local Dirichlet form on
$L^2(F;\rho\cdot\mu)$ in the sense of~\cite{17}, Chapter \textup{IV},
Definition~3.1.
\end{thm}

\begin{pf}
The assertion follows from the main result in~\cite{27}.
\end{pf}

By virtue of Theorem~\ref{thm2.2} and~\cite{17}, there exists a~diffusion
process $M^\rho=(\Omega,\mathcal{M},\{\mathcal{M}_t\},\theta
_t,X_t,$ $P_z)$ on $F$ associated with the Dirichlet form $(\mathcal
{E}^\rho, \mathcal{F}^\rho)$. $M^\rho$ will be called distorted
OU-process on $F$. Since constant functions are in~$\mathcal{F}^\rho$
and $\mathcal{E}^\rho(1,1)=0$, $M^\rho$ is recurrent and
conservative. We denote by $\mathbf{A}_+^\rho$ the set of all
positive continuous additive functionals (PCAF in abbreviation) of
$M^\rho$, and define $\mathbf{A}^\rho:=\mathbf{A}^\rho_+-\mathbf
{A}^\rho_+$. For $A\in\mathbf{A}^\rho$, its total variation process
is denoted by $\{A\}$. We also define $\mathbf{A}^\rho_0:=\{A\in
\mathbf{A}^\rho|E_{\rho\cdotp\mu}(\{A\}_t)<\infty\ \forall t>0\}$.
Each element in $\mathbf{A}^\rho_+$ has a~corresponding positive
$\mathcal{E}^\rho$-smooth measure on $F$ by the Revuz correspondence.
The set of all such measures will be denoted by~$S^\rho_+$.
Accordingly, $A_t\in\mathbf{A}^\rho$ corresponds to a~$\nu\in
S^\rho:=S^\rho_+-S^\rho_+$, the set of all $\mathcal{E}^\rho
$-smooth signed measure in the sense that $A_t=A_t^1-A_t^2$ for\vspace*{1pt}
$A_t^k\in\mathbf{A}^\rho_+,k=1,2$, whose Revuz measures are $\nu^k,
k=1,2$, and $\nu=\nu^1-\nu^2$ is the Hahn--Jordan decomposition of
$\nu$. The element of $\mathbf{A}^\rho$ corresponding to $\nu\in
S^\rho$ will be denoted by$A^\nu$.

Note that for each $l\in H$ the function $u(z)=\langle l,z\rangle$
belongs to the extended Dirichlet space $\mathcal{F}^\rho_e$ and
%
\begin{equation} \label{eq2.1}
\mathcal{E}^\rho(l(\cdot),v)=\frac{1}{2}\int\langle l,
Dv(z)\rangle\rho(z)\,d\mu(z)\qquad \forall v\in C_b^1(F).
\end{equation}
On the other hand, the AF $\langle l,X_t-X_0\rangle$ of $M^\rho$
admits a~unique decomposition into a~sum of a~martingale AF ($M_t$) of
finite energy and CAF ($N_t$) of zero energy. More precisely, for every
$l\in H$,
%
\begin{equation} \label{eq2.2}
\langle l,X_t-X_0\rangle=M^l_t+N^l_t\qquad \forall t\geq0\
P_z\mbox{-a.s.}
\end{equation}
for $\mathcal{E}^\rho$-q.e. $z\in F$.

Now for $\rho\in L^1(H,\mu)$ and $l\in H$, we say that $\rho\in
\operatorname{BV}_l(H)$ if there exists a~constant $C_l>0$,
%
\begin{equation} \label{eq2.3}
\biggl|\int\langle l, Dv(z)\rangle\rho(z)\,d\mu(z)\biggr|\leq C_l\|
v\|_\infty\qquad \forall v\in C_b^1(F).
\end{equation}

By the same argument as in~\cite{13}, Theorem 2.1, we obtain the following
theorem.

\begin{thm}\label{thm2.3}
Let $\rho\in L^1_+$ and $l\in H$.
\begin{longlist}[(2)]
\item[(1)] The following two conditions are equivalent:
\begin{longlist}[(iii)]
\item[(i)]$\rho\in \operatorname{BV}_l(H)$.

\item[(ii)] There exists a~(unique) signed measure $\nu_l$ on $F$ of finite
total variation such that
%
\begin{equation} \label{eq2.4}
\frac{1}{2}\int\langle l,  Dv(z)\rangle\rho(z)\,d\mu(z)=-\int
_Fv(z)\nu_l(dz) \qquad \forall v\in C_b^1(F).
\end{equation}
In this case, $\nu_l$ necessarily belongs to $S^{\rho+1}$.

Suppose further that $\rho\in \operatorname{QR}(H)$. Then the following condition is
also equivalent to the above:

\item[(iii)]$N^l\in\mathbf{A}_0^\rho$.
\end{longlist}
In this case, $\nu_l\in S^\rho$ and $N^l=A^{\nu_l}$.

\item[(2)] $M^l$ is a~martingale AF with quadratic variation process
%
\begin{equation} \label{eq2.5}
\langle M^l\rangle_t=t|l|^2,\qquad t\geq0.
\end{equation}
\end{longlist}
\end{thm}

\begin{remark}\label{rem2.4}
Recall that the Riesz representation theorem of
positive linear functionals on continuous functions by measures is not
applicable to obtain Theorem~\ref{thm2.3}, (i)${}\Rightarrow{}$(ii), because of the
lack of local compactness. However, the quasi-regularity of the
Dirichlet form provides a~means to circumvent this difficulty.
\end{remark}

In the rest of this section, we shall introduce a~special class of
$\rho\in \operatorname{QR}(H)$, which will be used in Section~\ref{sec4} below.

A nonnegative measurable function $h(s)$ on $\mathbb{R}^1$ is said to
possess the \emph{Hamza property} if $h(s)=0$  $ds$-a.e. on
the closed set $\mathbb{R}^1\setminus R(h)$ where
\[
R(h)=\biggl\{s\in\mathbb{R}^1\dvtx \int_{s-\varepsilon}^{s+\varepsilon}\frac
{1}{h(r)}\,dr<\infty\mbox{ for some } \varepsilon>0\biggr\}.
\]
We say that a~function $\rho\in L^1_+(H,\mu)$ satisfies \emph{the
ray Hamza condition in direction $l\in H$} ($\rho\in\mathbf{H}_l$ in
notation) if there exists a~nonnegative function $\tilde{\rho}_l$
such that
\[
\tilde{\rho}_l=\rho\ \mu\mbox{-a.e.} \mbox{ and }\tilde{\rho
}_l(z+sl) \mbox{ has the Hamza property in } s\in\mathbb{R}^1
\mbox{ for each } z\in H .
\]
We set $\mathbf{H}:=\bigcap_{k} \mathbf{H}_{e_k}$, where $e_k$ is as in
Hypothesis~\ref{hypo2.1}. A function in the family $\mathbf{H}$ is simply said
to satisfy the\emph{ ray Hamza condition}. By~\cite{5} $\mathbf
{H}\subset \operatorname{QR}(H)$, and thus we always have $\rho+1\in \operatorname{QR}(H)$, since
clearly $\rho+1\in\mathbf{H}$.

Next, we will present some explicit description of the Dirichlet form
$(\mathcal{E}^\rho, \mathcal{F}^\rho)$ for $\rho\in\mathbf{H}$.

For $e_j\in H$ as in Hypothesis~\ref{hypo2.1}, we set $H_{e_j}=\{se_j\dvtx s\in
\mathbb{R}^1\}$. We then have the direct sum decomposition
$H=H_{e_j}\oplus E_{e_j}$ given by
\[
z=se_j+x,\qquad s= \langle e_j,z \rangle.\vadjust{\goodbreak}
\]
Let $\pi_j$ be the projection onto the space $E_{e_j}$ and\vspace*{-1pt} $\mu
_{e_j}$ be the image measure of~$\mu$ under $\pi_j\dvtx H\rightarrow
E_{e_j}$, that is, $\mu_{e_j}=\mu\circ\pi_j^{-1}$. Then we see that for
any $F\in L^1(H,\mu)$
%
\begin{equation} \label{eq2.6}
\int_HF(z)\mu(dz)=\int_{E_{e_j}}\int_{\mathbb
{R}^1}F(se_j+x)p_j(s)\,ds\,\mu_{e_j}(dx),
\end{equation}
where $p_j(s)=(1/\sqrt{2\pi\lambda_j})e^{-s^2/2\lambda_j}$.
Thus by ~\cite{5}, Theorem 3.10, for all $u,v\in D(\mathcal{E}^\rho)$,
%
\begin{equation} \label{eq2.7}
\mathcal{E}^\rho(u,v)=\sum_{j=1}^\infty\mathcal{E}^{\rho
,e_j}(u,v),
\end{equation}
where
%
\begin{eqnarray} \label{eq2.8}
\mathcal{E}^{\rho,e_j}(u,v)&=&\frac{1}{2}\int_{E_{e_j}}\int_{R(\rho
(\cdot e_j+x))}\frac{d\tilde{u}_j(se_j+x)}{ds}\times\frac{d\tilde
{v}_j(se_j+x)}{ds}\nonumber\\[-8pt]\\[-8pt]
&&\hphantom{\frac{1}{2}\int_{E_{e_j}}\int_{R(\rho(\cdot e_j+x))}}
{}\times\rho(se_j+x)p_j(s)\,ds\,\mu_{e_j}(dx),\nonumber
\end{eqnarray}
and $u,\tilde{u}_j$ satisfy $\tilde{u}_j=u $ $\rho\mu$-a.e.
and $\tilde{u}_j(se_j+x)$ is absolutely continuous in~$s$ on $R(\rho
(\cdotp e_j+x))$ for each $x\in E_{e_j}$. $v$ and $\tilde{v}_j$ are
related in the same way.

\section{BV functions and distorted OU-processes in $F$}\label{sec3}
As in~\cite{13}, we introduce some function spaces on~$H$. Let
\[
A_{1/2}(x):=\int_0^x\bigl(\log(1+s)\bigr)^{1/2}\,ds,\qquad x\geq0,
\]
and let $\psi$ be its complementary function, namely,
\[
\psi(y):=\int_0^y(A'_{1/2})^{-1}(t)\,dt=\int_0^y\bigl(\exp(t^2)-1\bigr)\,dt.
\]
Define
\begin{eqnarray*}
L(\log L)^{1/2}(H,\mu)&:=&\{f\dvtx H\rightarrow\mathbb{R} |f \mbox{
Borel measurable, }A_{1/2}(|f|)\in L^1(H,\mu)\},
\\
L^\psi(H,\mu)&:=&\{g\dvtx H\rightarrow\mathbb{R}|g \mbox{ Borel
measurable, }\psi(c|g|)\in L^1(H,\mu) \\
&&\hspace*{167pt}\mbox{for some } c>0\}.
\end{eqnarray*}
From the general theory of Orlicz spaces (cf.~\cite{24}), we have the
following properties:
\begin{longlist}[(iii)]
\item[(i)] $L(\log L)^{1/2}$ and $L^\psi$ are Banach spaces under the norms
\begin{eqnarray*}
\|f\|_{L(\log L)^{1/2}}&=&\inf\biggl\{\alpha>0\Big|\int_H A_{1/2}(|f|/\alpha)\,d\mu\leq1\biggr\},
\\
\|g\|_{L^\psi}&=&\inf\biggl\{\alpha>0\Big|\int_H\psi(|g|/\alpha)\,d\mu\leq1\biggr\}.
\end{eqnarray*}
\item[(ii)] For $f\in L(\log L)^{1/2}$ and $g\in L^\psi$, we have
%
\begin{equation} \label{eq3.1}
\|fg\|_1\leq2\|f\|_{L(\log L)^{1/2}}\|g\|_{L^\psi}.
\end{equation}
\item[(iii)] Since~$\mu$ is Gaussian, the function $x\mapsto\langle
x,l\rangle$ belongs to $L^\psi$.
\end{longlist}

Let $c_j ,j\in\mathbb{N}$, be a~sequence in $[1,\infty)$. Define
\[
H_1:=\Biggl\{x\in H\Big|\sum_{j=1}^\infty\langle x,e_j\rangle^2c_j^2<\infty\Biggr\},
\]
equipped with the inner product
\[
\langle x,y\rangle_{H_1}:=\sum_{j=1}^\infty c_j^2\langle x,e_j\rangle
\langle y,e_j\rangle.
\]
Then clearly $(H_1, \langle\cdot,\cdot\rangle_{H_1})$ is a~Hilbert space such
that $H_1\subset H$ continuously and densely. Identifying~$H$ with its
dual we obtain the continuous and dense embeddings
\[
H_1\subset H(\equiv H^*)\subset H^*_1.
\]
It follows that
\[
_{H_1}\langle z,v\rangle_{H_1^*}=\langle z,v\rangle_H \qquad \forall z\in
H_1,v\in H,
\]
and that
$(H_1,H,H_1^*)$ is a~Gelfand triple.\vspace*{-1pt} Furthermore, $\{\frac{e_j}{c_j}\}
$ and $\{c_je_j\}$ are orthonormal bases of $H_1$ and $H_1^*$, respectively.

We also introduce a~family of~$H$-valued functions on~$H$ by
\[
(C_b^1)_{D(A)\cap H_1}:=\Biggl\{G\dvtx G(z)=\sum_{j=1}^m g_j(z)l^j,z\in H, g_j\in
C_b^1(H),l^j\in D(A)\cap H_1\Biggr\}.
\]
Denote by $D^*$ the adjoint of $D\dvtx C_b^1(H)\subset L^2(H,\mu
)\rightarrow L^2(H,\mu;H)$. That is
\begin{eqnarray*}
\operatorname{Dom}(D^*)&:=&\biggl\{\!G\in L^2(H,\mu;H)\Big|\\
&&\hphantom{\biggl\{\!}C_b^1\ni u\!\mapsto\!\int\!\langle
G,Du\rangle \,d\mu\mbox{ is continuous with respect to } L^2(H,\mu)\!\biggr\}.
\end{eqnarray*}
Obviously, $(C_b^1)_{D(A)\cap H_1}\subset \operatorname{Dom}(D^*)$. Then
%
\begin{eqnarray}\label{eq3.2}
\int_H D^*G(z)f(z)\mu(dz)=\int_H\langle G(z), Df(z)\rangle\mu
(dz)\nonumber\\[-8pt]\\[-8pt]
\eqntext{\forall G\in(C_b^1)_{D(A)\cap H_1},f\in C_b^1(H).}
\end{eqnarray}

For $\rho\in L(\log L)^{1/2}(H,\mu)$, we set
\[
V(\rho):=\sup_{G\in(C_b^1)_{D(A)\cap H_1},\|G\|_{H_1}\leq1}\int_H
D^*G(z)\rho(z)\mu(dz).\vadjust{\goodbreak}
\]
A function $\rho$ on~$H$ is called a~BV function in the Gelfand triple
$(H_1, H, H^*_1)$ [$\rho\in \operatorname{BV}(H, H_1)$ in notation], if $\rho\in
L(\log L)^{1/2}(H,\mu)$ and $V(\rho)$ is finite.
When $H_1=H=H_1^*$, this coincides with the definition of BV functions
defined in~\cite{1} and clearly $\operatorname{BV}(H,H)\subset \operatorname{BV}(H,H_1)$. We can
prove the following theorem by a~modification of the proof of~\cite{12},
Theorem~\ref{thm3.1}.

\setcounter{thm}{-1}
\begin{remark}\label{rem3.0}
The introduction of BV functions in a~Gelfand triple
is natural and originates from standard ideas when working with
infinite dimensional state spaces. The intersection of $\operatorname{BV}_l(H)$, when
$l$ runs through $D(A)\cap H_1$, describes functions which are
``componentwise of bounded variation'' in the sense that their weak
partial derivatives are measures. In contrast to finite dimensions,
this does not give rise
to vector-valued measures representing their total weak derivatives or
gradients. Therefore, one introduces an appropriate ``tangent space''
$H_1^*$ to~$H$, in which these total derivatives can be
represented as a~$H_1^*$-valued measure. This approach substantially
extends the applicability of the theory of BV functions on Hilbert
spaces. We document this by including the well-studied case of
linear SPDE with reflection, more precisely, the randomly vibrating
Gaussian string, forced to stay above a~level $\alpha\geq0$ (see
\cite{19,28}), which (in the case of $\alpha>0$) is then just a~special case of our general approach.
\end{remark}

\begin{thm}\label{thm3.1}
\textup{(i)} $\operatorname{BV}(H, H_1)\subset\bigcap_{l\in D(A)\cap H_1} \operatorname{BV}_l(H)$.
{\smallskipamount=0pt
\begin{longlist}[(iii)]
\item[(ii)] Suppose $\rho\in \operatorname{BV}(H, H_1)\cap L^1_+(H,\mu)$, then there exist
a positive finite measure $\|d\rho\|$ on~$H$ and a~Borel-measurable
map $\sigma_\rho\dvtx H\rightarrow H_1^*$ such that $\|\sigma_\rho(z)\|
_{H_1^*}=1 \|d\rho\|$-a.e, $\|d\rho\|(H)=V(\rho)$,
%
\begin{eqnarray}\label{eq3.3}
\int_HD^*G(z)\rho(z)\mu(dz)=\int_H { }_{H_1} \langle G(z), \sigma
_\rho(z)\rangle_{H_1^*}\|\,d\rho\|(dz)\nonumber\\[-8pt]\\[-8pt]
\eqntext{\forall G\in(C_b^1)_{D(A)\cap H_1}}
\end{eqnarray}
and $\|d\rho\|\in S^{\rho+1}$.

Furthermore, if $\rho\in \operatorname{QR}(H)$, $\|d\rho\|$ is $\mathcal{E}^\rho
$-smooth in the sense that it charges no set of zero $\mathcal
{E}_1^\rho$-capacity. In particular, the domain of integration~$H$ on
both sides of~(\ref{eq3.3}) can be replaced by $F$, the topological support of
$\rho\mu$.

Also, $\sigma_\rho$ and $\|d\rho\|$ are uniquely determined, that
is, if there are $\sigma_\rho'$ and $\|d\rho\|'$ satisfying relation
(\ref{eq3.3}), then $\|d\rho\|=\|d\rho\|'$ and $\sigma_\rho(z)=\sigma_\rho
'(z)$ for $\|d\rho\|$-a.e. z.

\item[(iii)] Conversely, if equation~(\ref{eq3.3}) holds for $\rho\in L(\log L)^{1/2}(H,\mu
)$ and for some positive finite measure $\|d\rho\|$ and a~map $\sigma
_\rho$ with the stated properties, then $\rho\in \operatorname{BV}(H,H_1)$ and
$V(\rho)=\|d\rho\|(H)$.

\item[(iv)] Let $W^{1,1}(H)$ be the domain of the closure of $(D, C_b^1(H))$
with norm
\[
\|f\|:=\int_H \bigl(|f(z)|+|Df(z)|\bigr)\mu(dz).
\]
Then $W^{1,1}(H)\subset \operatorname{BV}(H,H)$ and equation~(\ref{eq3.3}) is satisfied for each
$\rho\in W^{1,1}(H)$. Furthermore,
\[
\|d\rho\|=|D\rho|\cdot\mu,\qquad  V(\rho)=\int_H|D\rho|\mu(dz),\qquad
\sigma_\rho=\frac{1}{|D\rho|}D\rho I_{\{|D\rho|>0\}}.
\]
\end{longlist}}
\end{thm}

\begin{pf}
(i) Let $\rho\in \operatorname{BV}(H,H_1)$ and $ l\in D(A)\cap H_1$. Take
$G\in(C_b^1)_{D(A)\cap H_1}$ of the type
%
\begin{equation}\label{eq3.4}
G(z)=g(z)l,\qquad z\in H,  g\in C_b^1(H) .
\end{equation}
By~(\ref{eq3.2})
\begin{eqnarray*}
\int_HD^*G(z)f(z)\mu(dz)&=&\int_H\langle G(z),Df(z)\rangle\mu(dz)\\
&=&-\int_H\langle l,Dg(z)\rangle f(z)\mu(dz)\\
&&{}+2\int_H\langle
Al,z\rangle g(z)f(z)\mu(dz) \qquad \forall f\in C_b^1(H);
\end{eqnarray*}
consequently,
%
\begin{equation}\label{eq3.5}
D^*G(z)=-\langle l, Dg(z)\rangle+2g(z)\langle Al,z\rangle.
\end{equation}
Accordingly,
%
\begin{eqnarray}\label{eq3.6}
\int_H\langle l, Dg(z)\rangle\rho(z)\mu(dz)&=&-\int_HD^*G(z)\rho
(z)\mu(dz)\nonumber\\[-8pt]\\[-8pt]
&&{}+2\int_H\langle Al, z\rangle g(z)\rho(z)\mu(dz).\nonumber
\end{eqnarray}
For any $g\in C_b^1(H)$, satisfying $\|g\|_\infty\leq1$, by~(\ref{eq3.1}) the
right-hand side is dominated by
\[
V(\rho)\|l\|_{H_1}+4\|\rho\|_{L(\log L)^{1/2}}\|\langle Al,\cdot
\rangle\|_{L^\psi}<\infty,
\]
hence, $\rho\in \operatorname{BV}_l(H)$.

(ii) Suppose $\rho\in L^1_+(H,\mu)\cap \operatorname{BV}(H,H_1)$. By (i) and
Theorem~\ref{thm2.3} for each $l\in D(A)\cap H_1$, there exists a~finite signed
measure $\nu_l$ on~$H$ for which equation~(\ref{eq2.4}) holds. Define
\[
D^A_l\rho(dz):=2\nu_l(dz)+2\langle Al,z \rangle\rho(z)\mu(dz).
\]
In view of~(\ref{eq3.6}), for any $G$ of type~(\ref{eq3.4}), we have
%
\begin{equation}\label{eq3.7}
\int_HD^*G(z)\rho(z)\mu(dz)=\int_Hg(z)D_l^A\rho(dz),
\end{equation}
which in turn implies
%
\begin{equation}\label{eq3.8}
V(D^A_l\rho)(H)=\sup_{g\in C_b^1(H),  \|g\|_\infty\leq1}\int_H
g(z)D_l^A\rho(dz)\leq V(\rho)\|l\|_{H_1},
\end{equation}
where $V(D^A_l\rho)$ denotes the total variation measure of the signed
measure~$D_l^A\rho$.

For the orthonormal basis $\{\frac{e_j}{c_j}\}$ of $H_1$, we set
%
\begin{eqnarray}\label{eq3.9}
\gamma_\rho^A&:=&\sum_{j=1}^\infty2^{-j}V(D_{
{e_j}/{c_j}}^A\rho),\nonumber\\[-8pt]\\[-8pt]
v_j(z)&:=&\frac{d D_{
{e_j}/{c_j}}^A\rho(z)}{d\gamma_\rho^A(z)},\qquad z\in H , j\in\mathbb{N}.\nonumber
\end{eqnarray}
$\gamma_\rho^A$ is a~positive finite measure with $\gamma_\rho
^A(H)\leq V(\rho)$ and $v_j$ is Borel-measurable. Since $D_{
{e_j}/{c_j}}^A\rho$ belongs to $S^{\rho+1}$, so does $\gamma_\rho^A$. Then for
%
\begin{equation}\label{eq3.10}
G_n:=\sum_{j=1}^n g_j\frac{e_j}{c_j}\in(C_b^1)_{D(A)\cap H_1},\qquad n\in
\mathbb{N},
\end{equation}
by~(\ref{eq3.7}) the following equation holds:
%
\begin{equation}\label{eq3.11}
\int_HD^*G_n(z)\rho(z)\mu(dz)=\sum_{j=1}^n\int_H
g_j(z)v_j(z)\gamma_\rho^A(dz).
\end{equation}
Since $|v_j(z)|\leq2^j $ $ \gamma_\rho^A$-a.e. and
$C^1_b(H)$ is dense in $L^1(H,\gamma_\rho^A)$, we can find
$v_{j,m}\in C_b^1(H)$ such that
\[
\lim_{m\rightarrow\infty}v_{j,m}=v_j   \qquad \gamma_\rho
^A\mbox{-a.e.}
\]
Substituting
%
\begin{equation}\label{eq3.12}
g_{j,m}(z):=\frac{v_{j,m}(z)}{\sqrt{\sum
_{k=1}^nv_{k,m}(z)^2+1/m}}
\end{equation}
for $g_j(z)$ in~(\ref{eq3.10}) and~(\ref{eq3.11}) we get a~bound
\[
\sum_{j=1}^n\int_H g_{j,m}(z)v_j(z)\gamma_\rho^A(dz)\leq V(\rho),
\]
because $\|G_n(z)\|_{H_1}^2=\sum_{j=1}^n g_{j,m}(z)^2\leq1 $ $\forall z\in H$. By letting $m\rightarrow\infty$, we obtain
\[
\int_H\sqrt{\sum_{j=1}^nv_j(z)^2}\gamma_\rho^A(dz)\leq V(\rho
) \qquad \forall n\in\mathbb{N}.
\]
Now we define
%
\begin{equation}\label{eq3.13}
\| d\rho\|:=\sqrt{\sum_{j=1}^\infty v_j(z)^2}\gamma_\rho
^A(dz)
\end{equation}
and $\sigma_\rho\dvtx H\rightarrow H_1^*$ by
%
\begin{equation}\label{eq3.14}\qquad
\sigma_\rho(z)= \cases{
\displaystyle \sum_{j=1}^\infty\frac{v_j(z)}{\sqrt{\sum
_{k=1}^\infty v_k(z)^2}}\cdot c_je_j,&\quad if $\displaystyle z\in\Biggl\{\sum
_{k=1}^\infty v_k(z)^2>0\Biggr\}$,\vspace*{2pt}\cr
0,&\quad otherwise.
}
\end{equation}
Then
%
\begin{equation}\label{eq3.15}
\|d\rho\|(H)\leq V(\rho),\qquad \|\sigma_\rho(z)\|_{H_1^*}=1
\qquad \|d\rho\|\mbox{-a.e.},
\end{equation}
$\|d\rho\|$ is $S^{\rho+1}$-smooth and $\sigma_\rho$ is
Borel-measurable. By~(\ref{eq3.11}) we see that the desired equation~(\ref{eq3.3})
holds for $G=G_n$ as in~(\ref{eq3.10}).
It remains to prove~(\ref{eq3.3}) for any $G$ of type~(\ref{eq3.4}), that is, $G=g\cdot
l,g\in C_b^1(H),l\in D(A)\cap H_1$. In view of~(\ref{eq3.6}), equation~(\ref{eq3.3}) then reads
%
\begin{eqnarray}\label{eq3.16}
&&-\int_H\langle l,  Dg(z)\rangle\rho(z)\mu(dz)+2\int_Hg(z)\langle
Al,z\rangle\rho(z)\mu(dz)\nonumber\\[-8pt]\\[-8pt]
&&\qquad =\int_Hg(z)_{H_1}\langle l, \sigma_\rho(z)\rangle_{H_1^*}\|\,d\rho\|(dz).\nonumber
\end{eqnarray}
We set
\[
k_n:=\sum_{j=1}^n\langle l,e_j\rangle e_j=\sum_{j=1}^n\biggl\langle l,\frac
{e_j}{c_j}\biggr\rangle_{H_1} \frac{e_j}{c_j},\qquad  G_n(z):=g(z)k_n.
\]
Thus, $k_n\rightarrow l$ in $H_1$ and $Ak_n\rightarrow Al$ in~$H$ as
$n\rightarrow\infty$.
But then also
\[
\lim_{n\rightarrow\infty}\int_H \langle Dg,k_n\rangle\rho\, d\mu
=\int_H \langle Dg,l\rangle\rho\, d\mu
\]
and
\[
\biggl|\int_Hg(z)\langle Ak_n,z\rangle\rho(z)\mu(dz)-\int_Hg(z)\langle
Al,z\rangle\rho(z)\mu(dz)\biggr|
\]
\[
\leq2\|g\|_\infty\|\rho\|_{L(\log L)^{1/2}}\|\langle Ak_n-Al,\cdotp
\rangle\|_{L^\psi}.
\]

Furthermore,
\[
\lim_{n\rightarrow\infty}\int_Hg(z)_{H_1}\langle k_n,\sigma_\rho
(z)\rangle_{H_1^*}\|d\rho\|(dz)=\int_Hg(z)_{H_1}\langle l,\sigma
_\rho(z)\rangle_{H_1^*}\|d\rho\|(dz).
\]
So letting $n\rightarrow\infty$ yields~(\ref{eq3.16}).

If $\rho\in \operatorname{QR}(H)$, we can get the claimed result by the same
arguments as above.

Uniqueness follows by the same argument as~\cite{13}, Theorem 3.9.

(iii) Suppose $\rho\in L (\log)^{1/2}(H,\mu)$ and that equation~(\ref{eq3.3})
holds for some positive finite measure $\|d\rho\|$ and some map
$\sigma_\rho$ with the properties stated in~(ii). Then clearly
\[
V(\rho)\leq\|d\rho\|(H)
\]
and hence $\rho\in \operatorname{BV}(H,H_1)$. To obtain the converse inequality, set
\[
\sigma_j(z):=\langle c_je_j,\sigma_\rho(z)\rangle
_{H_1^*}=_{H_1}\biggl\langle\frac{e_j}{c_j},\sigma_\rho(z)\biggr\rangle
_{H_1^*},\qquad  j\in\mathbb{N}.
\]
Fix an arbitrary $n$. As in the proof of (ii), we can find functions
\[
v_{j,m}\in C_b^1(H),\qquad  \lim_{m\rightarrow\infty}v_{j,m}(z)=\sigma
_j(z) \qquad  \|d\rho\|\mbox{-a.e.}
\]
Define $g_{j,m}(z)$ by~(\ref{eq3.12}). Substituting $G_{n,m}(z):=\sum
_{j=1}^ng_{j,m}(z)\frac{e_j}{c_j}$ for $G(z)$ in~(\ref{eq3.3}) then yields
\[
\sum_{j=1}^n\int_H g_{j,m}(z)\sigma_j(z)\|d\rho\|(dz)\leq V(\rho).
\]
By letting $m\rightarrow\infty$, we get
\[
\int_H\sqrt{\sum_{j=1}^n\sigma_j(z)^2}\|d\rho\|(dz)\leq V(\rho
) \qquad  \forall n\in\mathbb{N}.
\]
We finally let $n\rightarrow\infty$ to obtain $\|d\rho\|(H)\leq
V(\rho)$.

(iv) Obviously the duality relation~(\ref{eq3.2}) extends to $\rho\in
W^{1,1}(H)$ replacing $f\in C_b^1(H)$. By defining $\|d\rho\|$ and
$\sigma_\rho(z)$ in the stated way, the extended relation~(\ref{eq3.2}) is
exactly~(\ref{eq3.3}).
\end{pf}

\begin{thm}\label{thm3.2}
Let $\rho\in \operatorname{QR}(H)\cap \operatorname{BV}(H,H_1)$ and consider the
measure $\|d\rho\|$ and $\sigma_\rho$ from Theorem~\ref{thm3.1}\textup{(ii)}. Then
there is an $\mathcal{E}^\rho$-exceptional set \mbox{$S\subset F$} such that
$\forall z\in F\setminus S$ under $P_z$ there exists an $\mathcal
{M}_t$-cylindrical Wiener process~$W^z$, such that the sample paths of
the associated distorted OU-process~$M^\rho$ on $F$ satisfy the
following: for $l\in D(A)\cap H_1$
%
\begin{eqnarray}\label{eq3.17}
\langle l,X_t-X_0\rangle&=&\int_0^t\langle l,dW_s^z\rangle+\frac
{1}{2}\int_0^t { }_{H_1} \langle l,\sigma_\rho(X_s)\rangle
_{H_1^*}\,dL_s^{\|d\rho\|}\nonumber\\[-8pt]\\[-8pt]
&&{}-\int_0^t\langle Al, X_s\rangle
\,ds\qquad \forall t\geq0 \ P_z\mbox{-a.s.}\nonumber
\end{eqnarray}
Here $L_t^{\|d\rho\|}$ is the real valued PCAF associated with $\|
d\rho\|$ by the Revuz correspondence.\vadjust{\goodbreak}

In particular, if $\rho\in \operatorname{BV}(H,H)$, then $\forall z\in F\setminus
S$, $l\in D(A)\cap H$
\begin{eqnarray}
\langle l,X_t-X_0\rangle=\int_0^t\langle l,dW_s^z\rangle+\frac
{1}{2}\int_0^t\langle l,\sigma_\rho(X_s)\rangle \,dL_s^{\|d\rho\|
}-\int_0^t \langle Al,X_s\rangle \,ds\nonumber  \\
\eqntext{\forall t\geq0\ P_z\mbox{-a.s.}}
\end{eqnarray}
\end{thm}

\begin{pf}
Let $\{e_j\}$ be the orthonormal basis of~$H$ introduced above.
Define for all $k\in\mathbb{N}$
%
\begin{eqnarray}\label{eq3.18}
W_k^z(t)&:=&\langle e_k,X_t-z\rangle-\frac{1}{2}\int_0^t { }_{H_1}
\langle e_k,\sigma_\rho(X_s)\rangle_{H_1^*}\,dL_s^{\|d\rho\|}\nonumber\\[-8pt]\\[-8pt]
&&{}+\int_0^t\langle Ae_k, X_s\rangle \,ds .\nonumber
\end{eqnarray}
By~(\ref{eq2.1}) and~(\ref{eq3.16}), we get for all $k\in\mathbb{N}$
\begin{eqnarray}
\mathcal{E}^\rho(e_k(\cdot),g)=\int_Hg(z)\langle Ae_k,z\rangle\rho
(z)\mu(dz)-\frac{1}{2}\int_Hg(z)_{H_1}\langle e_k,\sigma_\rho
(z)\rangle_{H_1^*}\|d\rho\|(dz)  \nonumber \\
\eqntext{\forall g\in C_b^1(H).}
\end{eqnarray}
By Theorem~\ref{thm2.3}, it follows that for all $k\in\mathbb{N}$
%
\begin{equation}\label{eq3.19}
N_t^{e_k}=\frac{1}{2}\int_0^t { }_{H_1} \langle e_k,\sigma_\rho
(X_s)\rangle_{H_1^*}\, dL_s^{\|d\rho\|}-\int_0^t\langle Ae_k,
X_s\rangle \,ds .
\end{equation}
Here we get from~(\ref{eq3.18}),~(\ref{eq3.19}) and the uniqueness of decomposition
(\ref{eq2.2}) that for $\mathcal{E}^\rho$-q.e. $ z\in F$,
\[
W^z_k(t)=M_t^{e_k} \qquad  \forall t\geq0 \ P_z\mbox{-a.s.}
,
\]
where the $\mathcal{E}^\rho$-exceptional set and the zero measure set
does not depend on~$e_k$.
Indeed, we can choose the capacity zero set $S=\bigcup_{j=1}^\infty S_j$,
where $S_j$ is the $\mathcal{E}^\rho$-exceptional set for $e_j$, and
for $z\in F\setminus S$, we can use the same method to get a~zero
measure set independent of $e_k$.
By Dirichlet form theory, we get $\langle M^{e_i},M^{e_j}\rangle
_t=t\delta_{ij}$. So for $z\in F\setminus S$, $W^z_k$ is an $\mathcal
{M}_t$-Wiener process under $P_z$. Thus, with $W^z$ being an $\mathcal
{M}_t$-cylindrical Wiener process given by $W^z(t)=(W_k^z(t)e_k)_{k\in
\mathbb{N}}$,~(\ref{eq3.17}) is satisfied for $P_z$-a.e., where $z\in
F\setminus S$.
\end{pf}

\section{Reflected OU-processes}\label{sec4}
In this section, we consider the situation where $\rho=I_\Gamma\in
\operatorname{BV}(H,H_1)$, where $\Gamma\subset H$ and
\[
I_\Gamma(x)= \cases{
1,&\quad if $x\in\Gamma$, \cr
0,&\quad if $x\in\Gamma^c $.
}
\]
Denote the corresponding objects $\sigma_\rho,\|dI_\Gamma\|$ in
Theorem~\ref{thm3.1}(ii) by $-\mathbf{n}_\Gamma,\|\partial\Gamma\|$,
respectively. Then formula~(\ref{eq3.3}) reads
\[
\int_\Gamma D^*G(z)\mu(dz)=-\int_F { }_{H_1} \langle G(z),\mathbf
{n}_\Gamma\rangle_{H_1^*}\|\partial\Gamma\|(dz) \qquad \forall
G\in(C_b^1)_{D(A)\cap H_1},
\]
where the domain of integration $F$ on the right-hand side is the
topological support of $I_\Gamma\cdot\mu$. $F$ is contained in $\bar
{\Gamma}$, but we shall show that the domain of integration on the
right-hand side can be restricted to $\partial\Gamma$. We need to use
the associated distorted OU-process $M^{I_\Gamma}$ on $F$, which will
be called reflected OU-process on $\Gamma$.

First, we consider a~$\mu$-measurable set $\Gamma\subset H$ satisfying
%
\begin{equation}\label{eq4.1}
I_\Gamma\in \operatorname{BV}(H,H_1)\cap\mathbf{H}.
\end{equation}

\begin{remark}\label{rem4.1}
We emphasize that if $\Gamma$ is a~convex closed set
in~$H$, then obviously $I_\Gamma\in\bf{H}$. Indeed, for each $z,l\in
H$ the set $\{ s\in\mathbb{R}| z+sl\in\Gamma\}$ is a~closed
interval in $\mathbb{R}$, whose indicator function hence trivially has
the Hamza property. Hence, in particular, $I_\Gamma\in \operatorname{QR}(H)$.
\end{remark}

By a~modification of~\cite{12}, Theorem 4.2, we can prove the following theorem.

\begin{thm}\label{thm4.2}
Let $\Gamma\subset H$ be~$\mu$-measurable
satisfying condition~(\ref{eq4.1}). Then the support of $\|\partial\Gamma\|$
is contained in the boundary $\partial\Gamma$ of $\Gamma$, and the
following generalized Gauss formula holds:
%
\begin{eqnarray}\label{eq4.2}
\int_\Gamma D^*G(z)\mu(dz)=-\int_{\partial\Gamma} { }_{H_1}
\langle G(z),\mathbf{n}_\Gamma\rangle_{H_1^*}\|\partial\Gamma\|
(dz)\nonumber\\[-8pt]\\[-8pt]
\eqntext{\forall G\in(C_b^1)_{D(A)\cap H_1}.}
\end{eqnarray}
\end{thm}

\begin{pf}
For any $G$ of type~(\ref{eq3.4}), we have from~(\ref{eq2.1}),~(\ref{eq3.5}) and~(\ref{eq3.7}) that
%
\begin{equation}\label{eq4.3}
\mathcal{E}^{I_\Gamma}(l(\cdot),g)-\int_\Gamma g(z)\langle
Al,z\rangle\mu(dz)=-\frac{1}{2}\int_Fg(z)D_l^AI_\Gamma(dz).
\end{equation}
Since the finite signed measure $D_l^AI_\Gamma$ charges no set of zero
$\mathcal{E}^{I_\Gamma}_1$-capacity, equation~(\ref{eq4.3}) readily extends to any
$\mathcal{E}^{I_\Gamma}$-quasicontinuous function $g\in\mathcal
{F}^{I_\Gamma}_b:=\mathcal{F}^{I_\Gamma}\cap L^\infty(\Gamma,\mu)$.

Denote by $\Gamma^0$ the interior of $\Gamma$. Then $\Gamma^0\subset
F\subset\bar{\Gamma}$. In view of the construction of the measure $\|
dI_\Gamma\|$ in Theorem~\ref{thm3.1}, it suffices to show that for $\frac
{e_j}{c_j}\in D(A)\cap H_1$
\[
V(D_{{e_j}/{c_j}}^AI_\Gamma)(\Gamma^0)=0.
\]
By linearity and since positive constants interchange with $\sup$, it
suffices to show that
%
\begin{equation}\label{eq4.4}
V(D_{e_j}^AI_\Gamma)(\Gamma^0)=0.
\end{equation}
Take an arbitrary $\varepsilon>0$ and set
\[
U:=\{z\in H\dvtx d(z,H\setminus\Gamma^0)>\varepsilon\},\qquad V:=\{z\in
H\dvtx d(z,H\setminus\Gamma^0)\geq\varepsilon\},
\]
where $d$ is the metric distance of the Hilbert space~$H$. Then $\bar
{U}\subset V$ and~$V$ is a~closed set contained in the open set $\Gamma
^0$. We define a~function $h$ by
%
\begin{equation}\label{eq4.5}
h(z):=1-E_z(e^{-\tau_V}),\qquad z\in F,
\end{equation}
where $\tau_V$ denotes the first exit time of $M^{I_\Gamma}$ from the
set $V$. The nonnegative function $h$ is in the space $\mathcal
{F}^{I_\Gamma}_b$ and furthermore it is $\mathcal{E}^{I_\Gamma
}$-quasicontinuous because it is $M^{I_\Gamma}$ finely continuous.

Moreover,
%
\begin{equation}\label{eq4.6}
h(z)>0\qquad \forall z\in U,\qquad  h(z)=0\qquad \forall z\in
F\setminus V.
\end{equation}
Set
%
\begin{equation}\label{eq4.7}
\nu_j(dz):=h(z)D^A_{e_j}I_\Gamma(dz)
\end{equation}
and
%
\begin{equation}\label{eq4.8}
I_g^j:=\mathcal{E}^{I_\Gamma}(e_j(\cdot),gh)-\int_\Gamma
g(z)h(z)\langle Ae_j,z\rangle\mu(dz).
\end{equation}
Then equation~(\ref{eq4.3}) with the $\mathcal{E}^{I_\Gamma}$-quasicontinuous
function $gh\in\mathcal{F}^{I_\Gamma}_b$ replacing $g$ implies
\[
I_g^j=-\frac{1}{2}\int_F g(z)\nu_j(dz).
\]
In order to prove~(\ref{eq4.4}), it is enough to show that $I_g^j=0$ for any
function $g(z)$ of the type
%
\begin{equation}\label{eq4.9}
\qquad g(z)=f(\langle e_j,z\rangle,\langle l_2,z\rangle,\ldots,\langle
l_m,z\rangle);\qquad l_2,\ldots,l_m\in H, f\in C_0^1(R^m)
\end{equation}
for we have then $\nu_j=0$.

On account of~(\ref{eq2.8}), we have the expression
%
\begin{eqnarray}\label{eq4.10}
\mathcal{E}^{I_\Gamma}(e_j(\cdot),gh)&=&\mathcal{E}^{I_\Gamma,e_j}(e_j(\cdot),gh)\nonumber\\[-8pt]\\[-8pt]
&=&\frac{1}{2}\int_{E_{e_j}}\int_{R_x}\frac{d(g\tilde{h})(se_j+x)}{ds}p_j(s)\,ds\mu_{e_j}(dx),\nonumber
\end{eqnarray}
where $R_x=R(I_\Gamma(\cdot e_j+x)), F_x:=\{s\dvtx se_j+x\in F\}$ for $x\in
E_{e_j}$ and $\tilde{h}$ is a~$I_\Gamma\cdot\mu$-version of $h$
appearing in the description of~(\ref{eq2.8}). For $x\in E_{e_j}$, set
\[
V_x:=\{s\dvtx se_j+x\in V\},\qquad
\Gamma_x^0:=\{s\dvtx se_j+x\in\Gamma^0\}.
\]
We then have the inclusion $V_x\subset\Gamma_x^0\subset R_x\cap F_x$.
By~(\ref{eq4.6}), $h(se_j+x)=0$ for any $x
\in E_{e_j}$ and for any $s\in R_x\setminus V_x$. On the other hand,
there exists a~Borel set $N\subset E_{e_j}$ with $\mu_{e_j}(N)=0$ such
that for each $x\in E_{e_j}\setminus N$,
\[
h(se_j+x)=\tilde{h}(se_j+x)\qquad  ds\mbox{-a.e.}
\]
Here we set $h\equiv0$ on $H\setminus F$.
Since $\tilde{h}(\cdot e_j+x)$ is absolutely continuous in $s$, we can
conclude that
\[
\tilde{h}(se_j+x)=0 \qquad\forall x\in E_{e_j}\setminus N,
 \forall s\in R_x\setminus V_x.\vadjust{\goodbreak}
\]

Fix $x\in E_{e_j}\setminus N$ and let $I$ be any connected component
of the one dimensional open set $R_x$. Furthermore, for any function
$g$ of type~(\ref{eq4.9}) we denote the support of $g(\cdot e_j+x)$ by $K_x$
(which is a~compact set) and choose a~bounded open interval $J$
containing $K_x$. Then $I\cap V_x\cap K_x$ is a~closed set contained in
the bounded open interval $I\cap J$ and
\[
g\tilde{h}(se_j+x)=0 \qquad \forall s\in(I\cap J)\setminus(I\cap
V_x\cap K_x).
\]
Therefore, an integration by part gives
\[
\int_{I\cap J}\frac{d(g\tilde{h})(se_j+x)}{ds}p_j(s)\,ds=\int_{I\cap
J}\frac{1}{\lambda_j}(g\tilde{h})(se_j+x)sp_j(s)\,ds.
\]
Combining this with~(\ref{eq4.8}) and~(\ref{eq4.10}), we arrive at
\begin{eqnarray*}
I_g^j&=&\int_{E_{e_j}}\int_{R_x}\frac{1}{2\lambda_j}(g\tilde{h})(se_j+x)sp_j(s)\,ds\,\mu_{e_j}(dx)\\
&&{}-\int_H g(z)h(z)\langle Ae_j,z\rangle I_\Gamma(z)\mu(dz)=0.
\end{eqnarray*}
\upqed
\end{pf}

Now we state Theorem~\ref{thm3.2} for $\rho=I_\Gamma$.

\begin{thm}\label{thm4.3}
Suppose $\Gamma\subset H$ is a~$\mu$-measurable set
satisfying condition~(\ref{eq4.1}). Then there is an $\mathcal{E}^\rho
$-exceptional set $S\subset F$ such that $\forall z\in F\setminus S$,
under $P_z$ there exists an $\mathcal{M}_t$-cylindrical Wiener
process $W^z$, such that the sample paths of the associated reflected
OU-process $M^\rho$ on $F$ with $\rho=I_\Gamma$ satisfy the
following: for $l\in D(A)\cap H_1$
%
\begin{eqnarray}\label{eq4.11}
\langle l,X_t-X_0\rangle&=&\int_0^t\langle l,dW_s^z\rangle-\frac
{1}{2}\int_0^t{ }_{H_1} \langle l,\mathbf{n}_\Gamma(X_s)\rangle
_{H_1^*}\,dL_s^{\|\partial\Gamma\|}\nonumber\\[-8pt]\\[-8pt]
&&{}-\int_0^t\langle Al, X_s\rangle
\,ds\qquad   P_z\mbox{-a.s.}\nonumber
\end{eqnarray}
Here, $L_t^{\|\partial\Gamma\|}$ is the real valued PCAF associated
with $\|\partial\Gamma\|$ by the Revuz correspondence, which has the
following additional property: $\forall z\in F\setminus S$
%
\begin{equation}\label{eq4.12}
I_{\partial\Gamma}(X_s)\,dL_s^{\|\partial\Gamma\|}=dL_s^{\|\partial
\Gamma\|}\qquad P_z\mbox{-a.s.}
\end{equation}

In particular, if $\rho\in \operatorname{BV}(H,H)$, then $\forall z\in F\setminus S,
l\in D(A)\cap H$
\begin{eqnarray*}
\langle l,X_t-X_0\rangle&=&\int_0^t\langle l,dW_s^z\rangle-\frac
{1}{2}\int_0^t\langle l,\mathbf{n}_\Gamma(X_s)\rangle \,dL_s^{\|
\partial\Gamma\|}\\
&&{}-\int_0^t \langle Al,X_s\rangle \,ds \qquad  \forall t\geq0\  P_z\mbox{-a.s.}
\end{eqnarray*}
\end{thm}

\begin{pf}
All assertions except for~(\ref{eq4.12}) follow from Theorem~\ref{thm3.2} for
$\rho:=I_\Gamma$. Equation~(\ref{eq4.12}) follows by Theorem~\ref{thm4.2} and~\cite{10}, Theorem
5.1.3.
\end{pf}

\section{Stochastic reflection problem on a~regular convex set}\label{sec5}
In this section, we consider $\Gamma$ satisfying~\cite{6}, Hypothesis
1.1(ii), with $K:=\Gamma$, that is:

\begin{hypo}\label{hypo5.1}
There exists a~convex $C^\infty$ function\vspace*{1pt} $g\dvtx
H\rightarrow\mathbb{R}$ with $g(0)=0$, $g'(0)=0$, and $D^2g$ strictly
positive definite, that is, $\langle D^2g(x)h$, $h \rangle\geq\gamma
|h|^2$ $\forall h\in H$ for some $\gamma>0$, such that
\[
\Gamma=\{x\in H\dvtx g(x)\leq1\},\qquad \partial\Gamma=\{x\in H\dvtx
g(x)=1\}.
\]
Moreover, we also suppose that $D^2g$ is bounded on $\Gamma$ and
$|Q^{1/2}Dg|^{-1}\in\bigcap_{p>1} L^p(H,\mu)$.
\end{hypo}

\begin{remark}\label{rem5.2}
By~\cite{6}, Lemma 1.2, $\Gamma$ is convex and closed and
there exists some constant $\delta>0$ such that $|Dg(x)|\leq\delta
\ \forall x\in\Gamma$.
\end{remark}

\subsection{Reflected OU processes on regular convex sets}\label{sec5.1}
Under Hypothesis~\ref{hypo5.1}, by~\cite{7}, Lemma A.1, we can prove that $I_\Gamma\in
\operatorname{BV}(H,H)\cap \operatorname{QR}(H)$:

\begin{thm}\label{thm5.3}
Assume that Hypothesis~\ref{hypo5.1} holds. Then $I_\Gamma\in
\operatorname{BV}(H,H)\cap \operatorname{QR}(H)$.
\end{thm}

\begin{pf} We first note that trivially by Remark~\ref{rem4.1} we have that
$I_\Gamma\in \operatorname{QR}(H)$. Let
\[
\rho_\varepsilon(x):=\exp\biggl(-\frac{(g(x)-1)^2}{\varepsilon} 1_{\{
g\geq1\}} \biggr),\qquad  x\in H.
\]
Thus,
\[
\lim_{\varepsilon\rightarrow0}\rho_\varepsilon=I_\Gamma.
\]
Moreover,
\[
D\rho_\varepsilon=-\frac{2}{\varepsilon}\rho_\varepsilon1_{\{
g\geq1\}}Dg(g-1) \qquad  \mu\mbox{-a.e.}
\]
By~\cite{7}, Lemma A.1, we have for $\varphi\in C_b^1(H)$
\begin{eqnarray*}
&&\lim_{\varepsilon\rightarrow0}\frac{1}{\varepsilon}\int_H\varphi
(x)1_{\{g(x)\geq1\}}\bigl(g(x)-1\bigr)\langle Dg(x),z \rangle\rho_\varepsilon
(x)\mu(dx)\\
&&\qquad =\frac{1}{2}\int_{\partial\Gamma} \varphi(y)\langle
n(y),z\rangle\frac{|Dg(y)|}{|Q^{1/2}Dg(y)|}\mu_{\partial\Gamma}(dy),
\end{eqnarray*}
where $n:=Dg/|Dg|$ is the exterior normal to $\partial\Gamma$ at $y$
and $\mu_{\partial\Gamma}$ is the surface measure on $\partial
\Gamma$ induced by~$\mu$ (cf.~\cite{6,7,16}),
whereas by~(\ref{eq3.2}) for any $\varphi\in C_b^1(H)$ and $z\in D(A)$
\begin{eqnarray*}
&&\lim_{\varepsilon\rightarrow0}\frac{1}{\varepsilon}\int
_H\varphi(x)1_{\{g(x)\geq1\}}\bigl(g(x)-1\bigr)\langle Dg(x),z \rangle\rho
_\varepsilon(x)\mu(dx)
\\
&&\qquad =-\lim_{\varepsilon\rightarrow0}\frac{1}{2}\int_H\langle D\rho
_\varepsilon(x),\varphi(x)z\rangle\mu(dx)
\\
&&\qquad =-\frac{1}{2}\lim_{\varepsilon\rightarrow0}\int_H \rho
_\varepsilon(x)D^*(\varphi z)(x)\mu(dx)
\\
&&\qquad =-\frac{1}{2}\int_H 1_\Gamma(x)D^*(\varphi z)(x)\mu(dx).
\end{eqnarray*}
Thus,
%
\begin{eqnarray}\label{eq5.1}
 &&\int_H 1_\Gamma(x)D^*(\varphi z)(x)\mu(dx)\nonumber\hspace*{-35pt}\\[-8pt]\\[-8pt]
 &&\qquad =-\int_{\partial\Gamma}
\varphi(x)\langle n(x),z\rangle\frac{|Dg(y)|}{|Q^{1/2}Dg(y)|}\mu
_{\partial\Gamma}(dx) \quad \forall z\in D(A),\varphi\in C_b^1.\nonumber\hspace*{-35pt}
\end{eqnarray}
By the proof of~\cite{7}, Lemma A.1, we get that $g$ is a~nondegenerate
map. So we can use the co-area formula (see~\cite{16}, Theorem 6.3.1,
Chapter V,
or~\cite{7}, (A.4)):
\[
\int_H f\mu(dx)=\int_0^\infty\biggl[\int_{g=r}f(y)\frac
{1}{|Q^{1/2}Dg(y)|}\mu_{\Sigma_r}(dy)\biggr]\,dr.
\]
By~\cite{16}, Theorem 6.2, Chapter  V, the surface measure is defined for all
$r\geq0$, moreover~\cite{16}, Theorem 1.1, Corollary 6.3.2, Chapter  V, imply
that $r\mapsto\mu_{\Sigma_r}$ is continuous in the topology induced
by $D_r^p(H)$ for some $p\in(1,\infty),r\in(0,\infty)$ (cf. \cite
{16}) on the measures on $(H, \mathcal{B}(H))$. Take $f\equiv1$ in
the co-area formula, then by the continuity property of the surface
measure with respect to $r$ we have that $\frac{1}{|Q^{1/2}Dg(y)|}\mu
_{\Sigma_r}(dy)$ is a~finite measure supported in $\{g=r\}$.\vspace*{-2pt} By Remark
\ref{rem5.2} and since $\mu_{\partial\Gamma}=\mu_{\Sigma_1}$, we have that\vspace*{1pt}
$ \frac{|Dg(y)|}{|Q^{1/2}Dg(y)|}\mu_{\partial\Gamma}$ is a~finite
measure. And hence by Theorem~\ref{thm3.1}(iii), we get $I_\Gamma\in \operatorname{BV}(H,H)$.
\end{pf}

Thus by Theorem~\ref{thm4.3}, we immediately get the following.

\begin{thm}\label{thm5.4}
Assume Hypothesis~\ref{hypo5.1}. Then there exists an $\mathcal
{E}^\rho$-exceptio\-nal set $S\subset F$ such that $\forall z\in
F\setminus S$, under $P_z$ there exists an $\mathcal{M}_t$-cylindrical Wiener process $W^z$, such that the sample paths of the
associated reflected OU-process $M^\rho$ on $F$ with $\rho=I_\Gamma$
satisfy the following: for $l\in D(A)\cap H_1$
\begin{eqnarray}
\langle l,X_t-X_0\rangle=\int_0^t\langle l,dW_s^z\rangle-\frac
{1}{2}\int_0^t\langle l,\mathbf{n}_\Gamma(X_s)\rangle \,dL_s^{\|
\partial\Gamma\|}-\int_0^t\langle Al, X_s\rangle\, ds\nonumber \\
\eqntext{\forall t\geq0  \ P_z\mbox{-a.e.,}}
\end{eqnarray}
where $\mathbf{n}_\Gamma:=\frac{Dg}{|Dg|}$ is the exterior normal to
$\Gamma$
and
\[
\|\partial\Gamma\|(dy)= \frac{|Dg(y)|}{|Q^{1/2}Dg(y)|}\mu_{\partial
\Gamma}(dy),
\]
where $\mu_{\partial\Gamma}$ is the surface measure induced by $\mu
$ (cf.~\cite{6,7,16}).
\end{thm}

\begin{remark}\label{rem5.5}
It can be shown that for $x\in\partial\Gamma$,
$\mathbf{n}_\Gamma(x)=\frac{Dg}{|Dg|}$ is the exterior normal to
$\Gamma$, that is, the unique element in~$H$ of unit length such that
\[
\langle\mathbf{n}_\Gamma(x),y-x \rangle\leq0 \qquad  \forall
y\in\Gamma.
\]
\end{remark}

\subsection{Existence and uniqueness of solutions}\label{sec5.2}
Let $\Gamma\subset H$ and our linear operator $A$ satisfy Hypotheses
\ref{hypo5.1} and~\ref{hypo2.1}, respectively. Consider the following stochastic
differential inclusion in the Hilbert space~$H$,
%
\begin{equation}\label{eq5.2}
\cases{
dX(t)+\bigl(AX(t)+N_\Gamma(X(t))\bigr)\,dt\ni dW(t),\cr
X(0)=x,
}
\end{equation}
where $W(t)$ is a~cylindrical Wiener process in~$H$ on a~filtered
probability space $(\Omega,\mathcal{F},\mathcal{F}_t,P)$ and
$N_\Gamma(x)$ is the normal cone to $\Gamma$ at $x$, that is,
\[
N_\Gamma(x)=\{z\in H\dvtx \langle z,y-x\rangle\leq0 \ \forall
y\in\Gamma\}.
\]

\begin{defin}\label{def5.6}
A pair of continuous $H\times\mathbb{R}$-valued
and $\mathcal{F}_t$-adapted processes $(X(t),L(t)),t\in[0,T]$, is
called a~solution of~(\ref{eq5.2}) if the following conditions hold:
\begin{longlist}[(ii)]
\item[(i)] $X(t)\in\Gamma$ for all $t \in[0,T]$ $P$-a.s.;

\item[(ii)] $L$ is an increasing process with the property that
\[
I_{\partial\Gamma}(X_s)\,dL_s=dL_s \qquad   P\mbox{-a.s.}
\]
and for any $l\in D(A)$ we have
\[
\langle l,X_t-x\rangle=\int_0^t\langle l,dW_s\rangle-\int
_0^t\langle l,\mathbf{n}_\Gamma(X_s)\rangle \,dL_s-\int_0^t\langle Al,
X_s\rangle\, ds\qquad  \forall t\geq0\  P\mbox{-a.s.,}
\]
where $\mathbf{n}_\Gamma$ is the exterior normal to $\Gamma$.
\end{longlist}
\end{defin}

\begin{remark}\label{rem5.7}
By Remark~\ref{rem5.5}, we know that $\mathbf{n}_\Gamma(x)\in
N_\Gamma(x)$ for all $x\in\partial\Gamma$. Hence by Definition~\ref{def5.6}(ii), it follows that Definition~\ref{def5.6} is appropriate to define a~solution for the multi-valued equation~(\ref{eq5.2}).
\end{remark}

We denote the semigroup with the infinitesimal generator $-A$ by
$S(t)$, $t\geq0$.

\begin{defin}\label{def5.8}
A pair of continuous $H\times\mathbb{R}$ valued
and $\mathcal{F}_t$-adapted processes $(X(t),L(t)),t\in[0,T]$, is
called a~mild solution of~(\ref{eq5.2}) if:
\begin{longlist}[(ii)]
\item[(i)] $X(t)\in\Gamma$ for all $t \in[0,T]$ $P$-a.s.;

\item[(ii)] $L$ is an increasing process with the property
\[
I_{\partial\Gamma}(X_s)\,dL_s=dL_s \qquad  P\mbox{-a.s.}
\]
and
\[
X_t=S(t)x+\int_0^tS(t-s)\,dW_s-\int_0^tS(t-s)\mathbf{n}_\Gamma
(X_s)\,dL_s \qquad \forall t\in[0,T] \ P\mbox{-a.s.,}
\]
where $\mathbf{n}_\Gamma$ is the exterior normal to $\Gamma$.
In particular, the appearing integrals have to be well defined.
\end{longlist}
\end{defin}

\begin{lemma}\label{lem5.9}
The process given by
\[
\int_0^tS(t-s)\mathbf{n}_\Gamma(X_s)\,dL_s
\]
is $P$-a.s. continuous and adapted to $\mathcal{F}_t, t\in[0,T]$.
This especially implies that it is predictable.
\end{lemma}

\begin{pf}
As $|S(t-s)\mathbf{n}_\Gamma(X_s)|\leq M_T|\mathbf{n}_\Gamma
(X_s)|, s\in[0,T]$, the integrals $\int_0^tS(t-s)\mathbf{n}_\Gamma
(X_s)\,dL_s, t\in[0,T],$ are well defined. For $0\leq s\leq t \leq T$,
\begin{eqnarray*}
&&\biggl|\int_0^sS(s-u)\mathbf{n}_\Gamma(X_u)\,dL_u-\int
_0^tS(t-u)\mathbf{n}_\Gamma(X_u)\,dL_u\biggr|\\
&&\qquad \leq\biggl|\int_0^s[S(s-u)-S(t-u)]\mathbf{n}_\Gamma(X_u)\,dL_u\biggr|+\biggl|\int
_s^tS(t-u)\mathbf{n}_\Gamma(X_u)\,dL_u\biggr|\\
&&\qquad \leq\int_0^s|[S(s-u)-S(t-u)]\mathbf{n}_\Gamma(X_u)|\,dL_u+\int
_s^t|S(t-u)\mathbf{n}_\Gamma(X_u)|\,dL_u,
\end{eqnarray*}
where the first summand converges to zero as $s\uparrow t$ or
$t\downarrow s$, because
\[
\bigl|1_{[0,s)}(u)[S(s-u)-S(t-u)]\mathbf{n}_\Gamma(X_u)\bigr|\rightarrow0
\qquad \mbox{as } s\uparrow t \mbox{ or } t\downarrow s.
\]
For the second summand, we have
\[
\int_s^t|S(t-u)\mathbf{n}_\Gamma(X_u)|\,dL_u\leq
M_T(L_t-L_s)\rightarrow0 \qquad \mbox{as } s\uparrow t \mbox{ or }
t\downarrow s.
\]
By the same arguments as in~\cite{23}, Lemma 5.1.9, we conclude that the
integral is adapted to $\mathcal{F}_t, t\in[0,T]$.
\end{pf}

\begin{thm}\label{thm5.10}
$(X(t),L_t), t\in[0,T]$, is a~solution of~(\ref{eq5.2}) if
and only if it is a~mild solution.
\end{thm}

\begin{pf}
$(\Rightarrow)$ First, we prove that for arbitrary $\zeta\in C^1([0,T],D(A))$ the
following equation holds:
%
\begin{eqnarray}\label{eq5.3}
\langle X_t,\zeta_t\rangle&=&\langle x,\zeta_0\rangle+\int
_0^t\langle\zeta_s,dW_s\rangle-\int_0^t\langle\mathbf{n}_\Gamma
(X_s),\zeta_s\rangle \,dL_s\nonumber\\[-8pt]\\[-8pt]
&&{}+\int_0^t\langle X_s,-A\zeta_s+\zeta'_s\rangle \,ds\qquad  \forall t\geq0\ P\mbox{-a.s.}\nonumber
\end{eqnarray}
If $\zeta_s=\eta f_s$ for $f\in C^1([0,T])$ and $\eta\in D(A)$,
by It\^{o}'s formula we have the above relation for such $\zeta$.
Then by~\cite{23}, Lemma G.0.10, and the same arguments as the proof of
Proposition G.0.11 we obtain the above formula for all $\zeta\in
C^1([0,T],D(A))$.
As in~\cite{23}, Proposition G.0.11, for the resolvent
$R_n:=(n+A)^{-1}\dvtx H\rightarrow D(A)$ and $t\in[0,T]$ choosing $\zeta
_s:=S(t-s)nR_n\eta, \eta\in H$, we deduce from~(\ref{eq5.3}) that
\begin{eqnarray}
\langle X_t,nR_n\eta\rangle&=&\langle x,S(t)nR_n\eta\rangle
+\int_0^t\langle S(t-s)nR_n\eta, dW_s\rangle\nonumber \\
&&{}-\int_0^t\langle
\mathbf{n}_\Gamma(X_s),S(t-s)nR_n\eta\rangle\, dL_s\nonumber \\
&&{}+\int_0^t\langle X_s, AS(t-s)nR_n\eta\rangle+\langle X_s, -AS(t-s)nR_n\eta
\rangle \,ds\nonumber \\
& =&\Biggl\langle S(t)x+\int_0^tS(t-s)\, dW_s+\int_0^t S(t-s)\mathbf{n}_\Gamma
(X_s) \,dL_s,nR_n\eta\Biggr\rangle\nonumber \\
\eqntext{\forall t\in[0,T]\
P\mbox{-a.s.}}
\end{eqnarray}
Letting $n\rightarrow\infty$, we conclude that $(X(t),L_t), t\in
[0,T],$ is a~mild solution.

$(\Leftarrow)$ By Lemma~\ref{lem5.9} and~\cite{23}, Theorem 5.1.3, we have
\[
\int_0^tS(t-s)\mathbf{n}_\Gamma(X_s)\,dL_s \quad \mbox{and}\quad  \int
_0^tS(t-s)\,dW_s,\qquad  t\in[0,T],
\]
have predictable versions. And we use the same notation for the
predictable versions of the respective processes.
As $(X_t,L_t)$ is a~mild solution, for all $\eta\in D(A)$ we get
\begin{eqnarray*}
\int_0^t\langle X_s, A\eta\rangle \,ds&=&\int_0^t\langle
S(s)x,A\eta\rangle \,ds\\
&&{}-\int_0^t\biggl\langle\int_0^sS(s-u)\mathbf
{n}_\Gamma(X_u)\,dL_u,A\eta\biggr\rangle \,ds\\
&&{}+\int_0^t\biggl\langle\int_0^sS(s-u)\,dW_u,A\eta\biggr\rangle \,ds\qquad \forall
t\in[0,T]\ P\mbox{-a.s.}
\end{eqnarray*}
The assertion that $(X(t),L_t), t\in[0,T],$ is a~solution of~(\ref{eq5.2}) now
follows as in the proof of~\cite{23}, Proposition G.0.9, because
\begin{eqnarray*}
&&\int_0^t\biggl\langle\int_0^sS(s-u)\mathbf{n}_\Gamma
(X_u)\,dL_u,A\eta\biggr\rangle \,ds\\[-3pt]
&&\qquad =\int_0^t\int_0^s\biggl\langle\mathbf
{n}_\Gamma(X_u),-\frac{d}{ds}S(s-u)\eta\biggr\rangle \,dL_u \,ds\\[-3pt]
&&\qquad =-\biggl\langle\int_0^t S(t-s)\mathbf{n}_\Gamma(X_s) \,dL_s,\eta\biggr\rangle
+\biggl\langle\int_0^t \mathbf{n}_\Gamma(X_s) \,dL_s,\eta\biggr\rangle.
\end{eqnarray*}
\upqed\vspace*{-3pt}
\end{pf}

Below, we prove~(\ref{eq5.2}) has a~unique solution in the sense of
Definition~\ref{def5.6}.\vspace*{-3pt}

\begin{thm}\label{thm5.11}
Let $\Gamma\subset H$ satisfy Hypothesis~\ref{hypo5.1}. Then
the stochastic inclusion~(\ref{eq5.2}) admits at most one solution in the sense
of Definition~\ref{def5.6}.\vspace*{-3pt}
\end{thm}

\begin{pf}
Let $(u,L^1)$ and $(v,L^2)$ be two solutions of~(\ref{eq5.2}), and let
$\{e_k\}_{k\in N}$ be the eigenbasis of $A$ from above. We then have
\begin{eqnarray*}
&&\langle e_k,u(t)-v(t)\rangle+\int_0^t\langle\alpha
_ke_k,u(s)-v(s)\rangle \,ds+\int_0^t\langle e_k,\mathbf{n}_\Gamma
(u(s))\rangle \,dL^1_s\\[-3pt]
&&\qquad {}-\int_0^t\langle e_k,\mathbf{n}_\Gamma
(v(s))\rangle \,dL^2_s=0.
\end{eqnarray*}
Setting $\phi_k(t):=\langle e_k,u(t)-v(t)\rangle$, we obtain
\begin{eqnarray}\label{eq5.4}
\phi_k^2(t)&=&2\int_0^t\phi_k(s)\,d\phi_k(s)\nonumber \\[-3pt]
&=&-2\biggl(\int_0^t\langle\alpha_ke_k,u(s)-v(s)\rangle\langle e_k,u(s)-v(s)\rangle
\,ds\nonumber \\[-3pt]
&&\hphantom{-2\biggl(}
{}+\int_0^t\langle e_k,\mathbf{n}_\Gamma(u(s))\rangle\langle
e_k,u(s)-v(s)\rangle \,dL^1_s\nonumber\\[-8.5pt]\\[-8.5pt]
&&\ \hphantom{-2\biggl(}{}-\int_0^t\langle e_k,\mathbf{n}_\Gamma
(v(s))\rangle\langle e_k,u(s)-v(s)\rangle \,dL^2_s\biggr)\nonumber \\[-3pt]
&\leq&-2\int
_0^t\langle e_k,\mathbf{n}_\Gamma(u(s))\rangle\langle
e_k,u(s)-v(s)\rangle \,dL^1_s\nonumber \\[-3pt]
&&{}+2\int_0^t\langle e_k,\mathbf{n}_\Gamma
(v(s))\rangle\langle e_k,u(s)-v(s)\rangle\, dL^2_s.\nonumber
\end{eqnarray}
By the dominated convergence theorem for all $t\geq0$, we have $P$-a.s.
\begin{eqnarray*}
&&\sum_{k\leq N}\int_0^t\langle e_k,\mathbf{n}_\Gamma(u(s))\rangle
\langle e_k,u(s)-v(s)\rangle \,dL^1_s
\\[-3pt]
&&\qquad \rightarrow\int_0^t\langle\mathbf{n}_\Gamma(u(s)),u(s)-v(s)\rangle
\,dL^1_s\qquad \mbox{as }N\rightarrow\infty
\end{eqnarray*}
and
\begin{eqnarray*}
&&\sum_{k\leq N}\int_0^t\langle e_k,\mathbf{n}_\Gamma(v(s))\rangle
\langle e_k,u(s)-v(s)\rangle\, dL^2_s
\\
&&\qquad \rightarrow\int_0^t\langle\mathbf{n}_\Gamma(v(s)),u(s)-v(s)\rangle
\,dL^2_s\qquad \mbox{as }N\rightarrow\infty.
\end{eqnarray*}
Summing over $k\leq N$ in~(\ref{eq5.4}) and letting $N\rightarrow\infty$
yield that for all $t\geq0$  $P$-a.s.
\begin{eqnarray*}
|u(t)-v(t)|^2&\leq&2\int_0^t\langle\mathbf{n}_\Gamma
(u(s)),v(s)-u(s)\rangle\, dL^1_s\\
&&{}+2\int_0^t\langle\mathbf{n}_\Gamma
(v(s)),u(s)-v(s)\rangle\, dL^2_s.
\end{eqnarray*}
By Remark~\ref{rem5.5} it follows that
\[
|u(t)-v(t)|^2\leq0,
\]
which implies
\[
u(t)=v(t),
\]
and thus
\[
\hspace*{142.5pt}L^1(t)=L^2(t).\hspace*{142.5pt}\qed
\]
\noqed
\end{pf}

Combining Theorems~\ref{thm5.4} and~\ref{thm5.11} with the Yamada--Watanabe theorem, we now
obtain the following theorem.

\begin{thm}\label{thm5.12}
If $\Gamma$ satisfies Hypothesis~\ref{hypo5.1}, then there
exists a~Borel set $M\subset H$ with $I_\Gamma\cdot\mu(M)=\mu
(\Gamma)$ such that for every $x\in M$,~(\ref{eq5.2}) has a~pathwise unique
continuous strong solution in the sense that for every probability
space $(\Omega,\mathcal{F},\mathcal{F}_t,P)$ with an $\mathcal
{F}_t$-Wiener process $W$, there exists a~unique pair of $\mathcal
{F}_t$-adapted processes $(X,L)$ satisfying Definition~\ref{def5.6} and
$P(X_0=x)=1$. Moreover, $X(t)\in M$ for all $t\geq0 $ $P$-a.s.
\end{thm}

\begin{pf}
By Theorems~\ref{thm5.4} and~\ref{thm5.11}, one sees that~\cite{15}, Theorem
3.14(a) is satisfied for the solution $(X,L)$. So, the assertion
follows from~\cite{15}, Theorem~3.14(b).
\end{pf}

\begin{remark}\label{rem5.13}
Following the same arguments as in the proof of~\cite{26},
Theorem 2.1, we can give an alternative proof of Theorem~\ref{thm5.12} for a~stronger notion of strong solutions (see, e.g.,~\cite{26}). Also,
because of Theorem~\ref{thm5.10}, by a~modification of~\cite{20}, Theorem 12.1, we
can prove the Yamada Watanabe theorem for the mild solution in
Definition~\ref{def5.8}, and then also a~corresponding version of Theorem~\ref{thm5.12}
for mild solutions for~(\ref{eq5.2}). This will be contained in forthcoming work.
\end{remark}

\subsection{The nonsymmetric case}\label{sec5.3}
In this section, we extend our results to the nonsymmetric case. For
$\Gamma\subset H$ satisfying Hypothesis~\ref{hypo5.1}, we consider the
nonsymmetric Dirichlet form,
\begin{eqnarray}
\mathcal{E}^\Gamma(u,v)=\int_\Gamma\biggl(\frac{1}{2}\langle
Du(z),Dv(z)\rangle+\langle B(z),Du(z)\rangle
v(z)\biggr)\mu(dz),\nonumber\\
\eqntext{u,v\in C_b^1(\Gamma),}
\end{eqnarray}
where $B$ is a~map from $\Gamma$ to~$H$ such that
%
\begin{eqnarray}\label{eq5.5}
&&B\in L^\infty(\Gamma\rightarrow H, \mu),\nonumber\\[-8pt]\\[-8pt]
&&\int_\Gamma\langle B,
Du\rangle \,d\mu\geq0 \qquad \mbox{for all } u\in C_b^1(\Gamma), u\geq0.\nonumber
\end{eqnarray}

Then $(\mathcal{E},C_b^1(\Gamma))$ is a~densely defined bilinear form
on $L^2(\Gamma;\mu)$ which is positive definite, since for all $u\in
C_b^1(\Gamma)$
\[
\mathcal{E}^\Gamma(u,u)=\int_\Gamma\frac{1}{2}\bigl(\langle
Du(z),Du(z)\rangle+\langle B(z),Du^2(z)\rangle(z)\bigr)\mu(dz)\geq0.
\]

Furthermore, by the same argument as~\cite{17}, Section II.3.e, we have $(\mathcal
{E},\allowbreak C_b^1(\Gamma))$ is closable on $L^2(\Gamma,\mu)$ and its closure
$(\mathcal{E}^\Gamma, \mathcal{F}^\Gamma)$ is a~Dirichlet form on
$L^2(\Gamma,\mu)$. We denote the extended Dirichlet space of
$(\mathcal{E}^\Gamma, \mathcal{F}^\Gamma)$ by $\mathcal
{F}_e^\Gamma$: Recall that $u\in\mathcal{F}_e^\Gamma$ if and only
if $|u|<\infty$ $I_\Gamma\cdot\mu$-a.e. and there exists a~sequence $\{u_n\}$
in~$\mathcal{F}^\Gamma$ such that $\mathcal
{E}^\Gamma(u_m-u_n,u_m-u_n)\rightarrow0$ as $n\geq m\rightarrow
\infty$ and $u_n\rightarrow u $ $I_\Gamma\cdot\mu$-a.e. as
$n\rightarrow\infty$.
This Dirichlet form satisfies the weak sector condition
\[
|\mathcal{E}_1^\Gamma(u, v)| \leq K\mathcal{E}_1^\Gamma(u,
u)^{1/2}\mathcal{E}_1^\Gamma(v, v)^{1/2}.
\]
Furthermore, we have the following theorem.

\begin{thm}\label{thm5.14}
Suppose $\Gamma\subset H$ satisfies Hypothesis~\ref{hypo5.1}.
Then $(\mathcal{E}^\Gamma, \mathcal{F}^\Gamma)$ is a~quasi-regular
local Dirichlet form on $L^2(\Gamma;\mu)$.
\end{thm}

\begin{pf}
The assertion follows by~\cite{17}, Section IV. 4b, and~\cite{28}.
\end{pf}

By virtue of Theorem~\ref{thm5.14} and~\cite{17}, there exists a~diffusion
process $M^\Gamma=(X_t,P_z)$ on $\Gamma$ associated with the
Dirichlet form $(\mathcal{E}^\Gamma, \mathcal{F}^\Gamma).$ Since
constant functions are in $\mathcal{F}^\Gamma$ and $\mathcal
{E}^\Gamma(1,1)=0$, $M^\Gamma$ is recurrent and conservative. We
denote by $\mathbf{A}_+^\Gamma$ the set of all positive continuous
additive functionals (PCAF in abbreviation) of $M^\Gamma$, and define
$\mathbf{A}^\Gamma=\mathbf{A}^\Gamma_+-\mathbf{A}^\Gamma_+$. For
$A\in\mathbf{A}^\Gamma$, its total variation process is denoted by
$\{A\}$. We also define $\mathbf{A}^\Gamma_0=\{A\in\mathbf
{A}^\Gamma|E_{I_\Gamma\cdotp\mu}(\{A\}_t)<\infty
\ \forall t>0\}$. Each element in $\mathbf{A}^\Gamma_+$ has a~corresponding positive $\mathcal{E}^\Gamma$-smooth measure on $\Gamma
$ by the Revuz correspondence. The totality of such measures will be
denoted by $S^\Gamma_+$. Accordingly, $\mathbf{A}^\Gamma$
corresponds to $S^\Gamma=S^\Gamma_+-S^\Gamma_+$, the set of
all
$\mathcal{E}^\Gamma$-smooth signed measure in the sense that
$A_t=A_t^1-A_t^2$ for\vadjust{\goodbreak} $A_t^k\in\mathbf{A}^\rho_+,k=1,2$, whose Revuz
measures are $\nu^k, k=1,2$, and $\nu=\nu^1-\nu^2$ is the
Hahn--Jordan decomposition of $\nu$. The element of $\mathbf{A}$
corresponding to $\nu\in S$ will be denoted by$A^\nu$.

Note that for each $l\in H$ the function $u(z)=\langle l,z\rangle$
belongs to the extended Dirichlet space $\mathcal{F}^\Gamma_e$ and
%
\begin{equation}\label{eq5.6}
\qquad \mathcal{E}^\Gamma(l(\cdot),v)=\int_\Gamma\biggl(\frac{1}{2}\langle l,
Dv(z)\rangle+\langle B(z),l\rangle v(z)\biggr)\mu(dz)\qquad  \forall
v\in C_b^1(\Gamma).
\end{equation}
On the other hand, the AF $\langle l,X_t-X_0\rangle$ of $M^\Gamma$
admits a~decomposition into a~sum of a~martingale AF $(M_t)$ of finite
energy and CAF $(N_t)$ of zero energy. More precisely, for every $l\in
H$
%
\begin{equation}\label{eq5.7}
\langle l,X_t-X_0\rangle=M^l_t+N^l_t \qquad  \forall t\geq0\
P_z\mbox{-a.s.}
\end{equation}
for $\mathcal{E}^\rho$-q.e. $z\in\Gamma$.

Then we have the following theorem.

\begin{thm}\label{thm5.15}
Suppose $\Gamma\subset H$ satisfies Hypothesis~\ref{hypo5.1}.
\begin{longlist}[(2)]
(1) The next three conditions are equivalent:
\begin{longlist}[(iii)]
\item[(i)]$N^l\in A_0$.

\item[(ii)]$|\mathcal{E}^\Gamma(l(\cdot),v)|\leq C\|v\|_\infty \ \forall v\in C_b^1(\Gamma)$.

\item[(iii)] There exists a~finite (unique) signed measure $\nu_l$ on $\Gamma
$ such that
%
\begin{equation}\label{eq5.8}
\mathcal{E}^\Gamma(l(\cdot),v)=-\int_\Gamma v(z)\nu_l(dz)\qquad  \forall v\in C_b^1(\Gamma).
\end{equation}
In this case, $\nu_l$ is automatically smooth and
\[
N^l=A^{\nu_l}.
\]
\end{longlist}

\item[(2)] $M^l$ is a~martingale AF with quadratic variation process
%
\begin{equation}\label{eq5.9}
\langle M^l\rangle_t=t|l|^2,\qquad t\geq0.
\end{equation}
\end{longlist}
\end{thm}

\begin{pf}
(1) By~\cite{21}, Theorem 5.2.7, and the same arguments as in~\cite{11}, we can extend Theorem 6.2 in~\cite{11} to our nonsymmetric case
to prove the assertions.

(2) Since
\[
\mathcal{E}^\Gamma(u,v)=\int_\Gamma\biggl(\frac{1}{2}\langle
Du(z),Dv(z)\rangle+\langle B(z),Du(z)\rangle v(z)\biggr)\mu(dz),\qquad u,v\in
\mathcal{F}^\Gamma,
\]
by~\cite{21}, Theorem 5.1.5, for $u\in C_b^1(\Gamma)$, $f\in\mathcal
{F}^\Gamma$ bounded we have
\begin{eqnarray*}
&&\int\tilde{f}(x)\mu_{\langle M^{[u]}\rangle
}(dx)\\
&&\qquad =2\mathcal{E}^\Gamma(u,uf)-\mathcal{E}^\Gamma(u^2,f)\\
&&\qquad =2\int_\Gamma\biggl(\frac{1}{2}\langle Du(z),D(u\tilde{f})(z)\rangle
+\langle B(z),Du(z)\rangle u(z)\tilde{f}(z)\biggr)\mu(dz)\\
&&\qquad \quad {}-\int_\Gamma
\biggl(\frac{1}{2}\langle D(u(z)^2),D\tilde{f}(z)\rangle+\langle
B(z),D(u^2)(z)\rangle\tilde{f}(z)\biggr)\mu(dz)\\
&&\qquad =
\int_\Gamma\langle Du(z),Du(z)\rangle\tilde{f}(z)\mu
(dz).
\end{eqnarray*}
Here $ \tilde{f}$ denotes the $\mathcal{E}^\Gamma$-quasi-continuous
version of $f$, $\mu_{\langle M^{[u]}\rangle}$ is the Reuvz measure
for $\langle M^{[u]}\rangle$ and $M^{[u]}$ is the martingale additive
functional in the Fukushima decomposition for $u(X_t)$.
Hence, we have
\[
\mu_{\langle M^{[u]}\rangle}(dz)=I_\Gamma\langle Du(z),Du(z)\rangle
\cdot\mu(dz).
\]

By~\cite{21}, (5.1.3), we also have
\[
e(\langle M^l\rangle)=e(M^l)
=\int_\Gamma\frac{1}{2}\langle l,l\rangle\mu(dz),
\]
where $e(M^l)$ is the energy of $M^l$. Then~(\ref{eq5.9}) easily
follows.
\end{pf}

By Theorem~\ref{thm3.1}, we can now prove the following theorem.

\begin{thm}\label{thm5.16}
Suppose $\Gamma\subset H$ satisfies Hypothesis~\ref{hypo5.1}.
Then there is an $\mathcal{E}^\Gamma$-exceptional set $S\subset
\Gamma$ such that $\forall z\in\Gamma\setminus S$, under $P_z$
there exists an $\mathcal{M}_t$-cylindrical Wiener process $W^z$,
such that the sample paths of the associated OU-process $M^\Gamma$ on
$\Gamma$ satisfy the following: for $l\in D(A)\cap H_1$
%
\begin{eqnarray}\label{eq5.11}
\langle l,X_t-X_0\rangle&=&\int_0^t\langle l,dW_s^z\rangle-\frac
{1}{2}\int_0^t{ }_{H_1} \langle\l,\mathbf{n}_\Gamma(X_s)\rangle
_{H_1^*\,}dL_s^{\|\partial\Gamma\|}\nonumber\\[-8pt]\\[-8pt]
&&{}-\int_0^t\langle Al, X_s\rangle
\,ds-\int_0^t\langle l,B(X_s)\rangle \,ds  \qquad P_z\mbox{-a.s.}\nonumber
\end{eqnarray}
Here, $L_t^{\|\partial\Gamma\|}$ is the real valued PCAF associated
with $\|\partial\Gamma\|$ by the Revuz correspondence, which has the
following additional property: $\forall z\in\Gamma\setminus S$
%
\begin{equation}\label{eq5.12}
I_{\partial\Gamma}(X_s)\,dL_s^{\|\partial\Gamma\|}=dL_s^{\|\partial
\Gamma\|}\qquad  P_z\mbox{-a.s.}
\end{equation}
Here $\mathbf{n}_\Gamma:=\frac{Dg}{|Dg|}$ is the exterior normal to
$\Gamma$,
and
\[
\|\partial\Gamma\|(dy)= \frac{|Dg(y)|}{|Q^{1/2}Dg(y)|}\mu_{\partial
\Gamma}(dy),
\]
where $\mu_{\partial\Gamma}$ the surface measure induced by~$\mu$.
\end{thm}

\begin{pf}
By~(\ref{eq5.6}) and~(\ref{eq3.16}), we have
\begin{eqnarray*}
\mathcal{E}^\Gamma(l(\cdot),v)&=&\int_\Gamma\frac{1}{2}\langle l,
Dv(z)\rangle+\langle B(z),l\rangle v(z)\mu(dz)
\\[-2pt]
&=&\int_\Gamma\langle B(z),l\rangle v(z)\mu(dz)+\int_\Gamma
v(z)\langle Al,z\rangle\mu(dz)\\[-2pt]
&&{}+\frac{1}{2}\int_{\partial\Gamma}
v(z)\langle l,\mathbf{n}_\Gamma(z)\rangle\|\partial\Gamma\|(dz).
\end{eqnarray*}
Thus, by Theorem~\ref{thm5.15}
\begin{eqnarray*}
N_t^l&=&-\biggl\langle Al,\int_0^t X_s(\omega)\,ds\biggr\rangle-\biggl\langle l,\int_0^t
B(X_s(\omega))\,ds\biggr\rangle\\[-2pt]
&&{}-\frac{1}{2}\biggl\langle l,\int_0^t\mathbf
{n}_\Gamma(X_s(\omega)) \,dL^{\|\partial\Gamma\|}_s(\omega)\biggr\rangle.
\end{eqnarray*}
By Theorem~\ref{thm5.15} and the same method as in Theorem~\ref{thm3.2} one then proves
the first assertion, and the last assertion follows by Theorems~\ref{thm5.3} and
\ref{thm5.4}.~%
\end{pf}

Let $\Gamma\subset H$ and our linear operator $A$ satisfy Hypotheses
\ref{hypo5.1} and~\ref{hypo2.1}, respectively. As in Section~\ref{sec5.2} we shall now
prove the existence and uniqueness of a~solution of the following
stochastic differential inclusion on the Hilbert space~$H$,
%
\begin{equation}\label{eq5.13}
\cases{
dX(t)+\bigl(AX(t)+B(X(t))+N_\Gamma(X(t))\bigr)\,dt\ni dW(t),\cr
X(0)=x,
}
\end{equation}
where $B$ satisfies condition~(\ref{eq5.5}), $W(t)$ is a~cylindrical Wiener
process in~$H$ on a~filtered probability space $(\Omega,\mathcal
{F},\mathcal{F}_t,P)$ and $N_\Gamma(x)$ is the normal cone to $\Gamma
$ at $x$, that is,
\[
N_\Gamma(x)=\{z\in H\dvtx \langle z,y-x\rangle\leq0 \ \forall
y\in\Gamma\}.
\]

\begin{defin}\label{def5.17}
A pair of continuous $H\times\mathbb{R}$-valued
and $\mathcal{F}_t$-adapted processes $(X(t),L(t))$, $t\in[0,T]$, is
called a~solution of~(\ref{eq5.13}) if the following conditions hold:
\begin{longlist}[(ii)]
\item[(i)] $X(t)\in\Gamma$ for all $t \in[0,T] $ $P$-a.s.;

\item[(ii)] $L$ is an increasing process with the property that
\[
I_{\partial\Gamma}(X_s)\,dL_s=dL_s\qquad  P\mbox{-a.s.},
\]
and for any $l\in D(A)$ we have
\begin{eqnarray*}
\langle l,X_t-x\rangle&=&\int_0^t\langle l,dW_s\rangle-\int
_0^t\langle l,\mathbf{n}_\Gamma(X_s)\rangle \,dL_s-\int_0^t \langle
l,B(X_s)\rangle\, ds\\[-2pt]
&&{}-\int_0^t\langle Al, X_s\rangle\, ds
\qquad \forall t\geq0 \ P\mbox{-a.s.},
\end{eqnarray*}
where $\mathbf{n}_\Gamma$ is the exterior normal to $\Gamma$.\vadjust{\goodbreak}
\end{longlist}

Below we prove~(\ref{eq5.13}) has a~unique solution in the sense of
Definition~\ref{def5.17}.\vspace*{-2pt}
\end{defin}

\begin{thm}\label{thm5.18}
Let $\Gamma\subset H$ satisfy Hypothesis~\ref{hypo5.1} and
$B$ satisfy the monotonicity condition
%
\begin{equation}\label{eq5.14}
\langle B(u)-B(v),u-v\rangle\geq-\alpha|u-v|^2
\end{equation}
for all $u,v\in \operatorname{dom}(B)$, for some $\alpha\in[0,\infty)$ independent
of $u,v$.
The stochastic inclusion~(\ref{eq5.13}) admits at most one solution in the
sense of Definition~\ref{def5.17}.\vspace*{-2pt}
\end{thm}

\begin{pf}
Let $(u,L^1)$ and $(v,L^2)$ be two solutions of~(\ref{eq5.13}), and let
$\{e_k\}_{k\in N}$ be the eigenbasis of $A$ from above. We then have
\begin{eqnarray*}
&&\langle e_k,u(t)-v(t)\rangle+\int_0^t\langle\alpha
_ke_k,u(s)-v(s)\rangle \,ds+\int_0^t\langle e_k,B(u(s))-B(v(s))\rangle \,ds
\\[-2pt]
&&\qquad {}+\int_0^t\langle e_k,\mathbf{n}_\Gamma(u(s))\rangle \,dL^1_s-\int
_0^t\langle e_k,\mathbf{n}_\Gamma(v(s))\rangle \,dL^2_s=0.
\end{eqnarray*}
Setting $\phi_k(t):=\langle e_k,u(t)-v(t)\rangle$, and we have
%
\begin{eqnarray}\label{eq5.15}
\phi_k^2(t)&=&2\int_0^t\phi_k(s)\,d\phi_k(s)\nonumber \\[-2pt]
&=&-2\biggl(\int_0^t\langle\alpha_ke_k,u(s)-v(s)\rangle\langle e_k,u(s)-v(s)\rangle
\,ds\\[-2pt]
&&\hphantom{-2\biggl(}
{}+\int_0^t\langle e_k,B(u(s))-B(v(s))\rangle\langle
e_k,u(s)-v(s)\rangle \,ds\nonumber\\[-2pt]
&&\hphantom{-2\biggl(}{}+\int_0^t\langle e_k,\mathbf{n}_\Gamma
(u(s))\rangle\langle e_k,u(s)-v(s)\rangle \,dL^1_s\nonumber\\[-2pt]
&&\hphantom{-2\biggl(}\hspace*{37pt}{}-\int_0^t\langle
e_k,\mathbf{n}_\Gamma(v(s))\rangle\langle e_k,u(s)-v(s)\rangle
\,dL^2_s\biggr)\nonumber \\[-2pt]
&\leq&-2\int_0^t\langle e_k,B(u(s))-B(v(s))\rangle\langle
e_k,u(s)-v(s)\rangle \,ds\nonumber \\[-2pt]
&&{}-2\int_0^t\langle e_k,\mathbf{n}_\Gamma
(u(s))\rangle\langle e_k,u(s)-v(s)\rangle \,dL^1_s\nonumber \\[-2pt]
&&{}+2\int_0^t\langle
e_k,\mathbf{n}_\Gamma(v(s))\rangle\langle e_k,u(s)-v(s)\rangle
\,dL^2_s.\nonumber
\end{eqnarray}
By the same argument as Theorem~\ref{thm5.11}, we have the following $P$-a.s.:
\begin{eqnarray*}
&&\sum_{k\leq N}\int_0^t\langle e_k,B(u(s))-B(v(s))\rangle\langle
e_k,u(s)-v(s)\rangle \,ds
\\[-2pt]
&&\qquad \rightarrow\int_0^t\langle B(u(s))-B(v(s)),u(s)-v(s)\rangle \,ds
\qquad \mbox{as }N\rightarrow\infty,
\\[-2pt]
&&\sum_{k\leq N}\int_0^t\langle e_k,\mathbf{n}_\Gamma(u(s))\rangle
\langle e_k,u(s)-v(s)\rangle \,dL^1_s
\\[-2pt]
&&\qquad \rightarrow\int_0^t\langle\mathbf{n}_\Gamma(u(s)),u(s)-v(s)\rangle
\,dL^1_s\qquad \mbox{as }N\rightarrow\infty
\end{eqnarray*}
and
\begin{eqnarray*}
&&\sum_{k\leq N}\int_0^t\langle e_k,\mathbf{n}_\Gamma(v(s))\rangle
\langle e_k,u(s)-v(s)\rangle \,dL^2_s
\\[-2pt]
&&\qquad \rightarrow\int_0^t\langle\mathbf{n}_\Gamma(v(s)),u(s)-v(s)\rangle
\,dL^2_s\qquad \mbox{as }N\rightarrow\infty.
\end{eqnarray*}
Summing over $k\leq N$ in~(\ref{eq5.15}) and letting $N\rightarrow\infty$
yield that for all $t\geq0$, $P$-a.s.
\begin{eqnarray*}
&&|u(t)-v(t)|^2+2\int_0^t\langle
B(u(s))-B(v(s)),u(s)-v(s)\rangle \,ds\\[-2pt]
&&\qquad \leq 2\int_0^t\langle\mathbf
{n}_\Gamma(u(s)),v(s)-u(s)\rangle \,dL^1_s+2\int_0^t\langle\mathbf
{n}_\Gamma(v(s)),u(s)-v(s)\rangle \,dL^2_s.
\end{eqnarray*}
By Remark~\ref{rem5.5}, it follows that
\[
|u(t)-v(t)|^2+2\int_0^t\langle B(u(s))-B(v(s)),u(s)-v(s)\rangle \,ds\leq0.
\]
By~(\ref{eq5.14}) and Gronwall's lemma, it follows that
\[
u(t)=v(t),
\]
and thus
\[
L^1(t)=L^2(t).
\]
\upqed
\end{pf}

Combining Theorems~\ref{thm5.16} and~\ref{thm5.18} with the Yamada--Watanabe theorem, we
obtain the following.

\begin{thm}\label{thm5.19}
If $\Gamma$ satisfies Hypothesis~\ref{hypo5.1} and $B$ in
(\ref{eq5.13}) satisfies~(\ref{eq5.14}), then there exists a~Borel set $M\subset H$
with $I_\Gamma\cdot\mu(M)=\mu(\Gamma)$ such that for every $x\in
M$,~(\ref{eq5.13}) has a~pathwise unique continuous strong solution in the
sense that for every probability space $(\Omega,\mathcal{F},\mathcal
{F}_t,P)$ with an $\mathcal{F}_t$-Wiener process $W$ there exists a~unique pair of $\mathcal{F}_t$-adapted processes $(X,L)$ satisfying
Definition~\ref{def5.17} and $P(X_0=x)=1$. Moreover, $X(t)\in M$ for all $t\geq
0 $ $P$-a.s.
\end{thm}

\begin{pf}
The proof is completely analogous to that of Theorem~\ref{thm5.12}.
\end{pf}

\section{Reflected OU-processeses on a~class of convex sets}\label{sec6}
Below for a~topological space $X$ we denote its Borel $\sigma$-algebra
by $\mathcal{B}(X)$. In this section, we consider\vadjust{\goodbreak} the case where
$H:=L^2(0,1)$, $\rho=I_{K_\alpha}$, where $K_\alpha:=\{f\in H|f\geq
-\alpha\},\alpha\geq0$, and $A=-\frac{1}{2}\frac{d^2}{dr^2}$ with
Dirichlet boundary conditions on $(0,1)$. So in this case $e_j=\sqrt
{2}\sin(j\pi r)$, $j\in\mathbb{N}$, is the corresponding
eigenbases. We recall that (cf.~\cite{28}) we have $\mu(C_0(
[0,1]))=1$. In~\cite{28}, L. Zambotti proved the following integration
by parts formulae in this situation:
\begin{itemize}
\item for $\alpha>0$,
\begin{eqnarray}
&&
\int_{K_\alpha}\langle l,D\varphi\rangle \,d\mu\nonumber\\[-2pt]
&&\qquad=-\int_{K_\alpha
}\varphi(x)\langle x,l''\rangle\mu(dx)-\int_0^1\,dr\,l(r)\int\varphi
(x)\sigma_\alpha(r,dx)\nonumber \\[-2pt]
\eqntext{\forall l\in D(A),\varphi\in C_b^1(H),}
\end{eqnarray}
\item for $\alpha=0$,
%
\begin{eqnarray}\label{eq6.1}
&&
\int_{K_0}\langle l,D\varphi\rangle\, d\nu\nonumber\\[-2pt]
&&\qquad=-\int_{K_0}\varphi
(x)\langle x,l''\rangle\nu(dx)-\int_0^1\,dr\,l(r)\int\varphi(x)\sigma
_0(r,dx)\\[-2pt]
\eqntext{\forall l\in D(A),\varphi\in C_b^1(H),}
\end{eqnarray}
\end{itemize}
where $\nu$ is the law of the Bessel Bridge of dimension 3 over $
[0,1]$ which is zero at $0$ and $1$, $\sigma_\alpha(r,dx)=\sigma
_\alpha(r)\mu_{\alpha}(r,dx),$ and for $\alpha>0$, $\sigma_\alpha
$ is a~positive bounded function, and for $\alpha=0$, $\sigma
_0(r)=\frac{1}{\sqrt{2\pi r^3(1-r)^3}},$ where $\mu_\alpha(r,dx),
\alpha\geq0,$ are probability kernels from $(H,\mathcal{B}(H))$ to
$([0,1],\mathcal{B}([0,1]))$.

\begin{remark}\label{rem6.1}
Since each $l$ in $D(A)$ has a~second derivative in $L^2$, its first
derivative is bounded, hence $l$ goes faster than linear to zero at any
point where $l$ is zero, in particular at the boundary points $r=0$ and $r=1$.
Hence, the second integral in the right-hand side of the above equality
is well defined.
\end{remark}

We know by~(\ref{eq3.5}) that for all $l\in D(A)$
\[
D^*(\varphi(\cdot)l)=-\langle l, D\varphi\rangle-\varphi\langle
l{''}, \cdot\rangle.
\]
Hence, for $\alpha>0$,
%
\begin{eqnarray}\label{eq6.2}
\int_{K_\alpha}D^*(\varphi(\cdot)l)\,d\mu=\int_0^1l(r)\int\varphi
(x)\sigma_\alpha(r,dx)\,dr \nonumber\\[-9pt]\\[-9pt]
\eqntext{\forall l\in D(A), \varphi\in C_b^1(H).}
\end{eqnarray}
Now take
%
\begin{equation}\label{eq6.3}
c_j:= \cases{
(j\pi)^{1/{2}+\varepsilon},&\quad if $\alpha>0$,\cr
(j\pi)^\beta,&\quad if $\alpha=0$,\vadjust{\goodbreak}
}
\end{equation}
where $\varepsilon\in(0,\frac{3}{2}]$ and $\beta\in(\frac
{3}{2},2]$, respectively,
and define
\[
H_1:=\Biggl\{x\in H\Big|\sum_{j=1}^\infty\langle x,e_j\rangle^2c_j^2<\infty\Biggr\}\vspace*{-2pt}
\]
equipped with the inner product
\[
\langle x,y\rangle_{H_1}:=\sum_{j=1}^\infty c_j^2\langle x,e_j\rangle
\langle y,e_j\rangle.\vspace*{-2pt}
\]
We note that $D(A)\subset H_1$ continuously for all $\alpha\geq0$,
since $\varepsilon\leq\frac{3}{2},\beta\leq2$.
Furthermore, $(H_1, \langle\cdot,\cdot\rangle_{H_1})$ is a~Hilbert space such
that $H_1\subset H$ continuously and densely. Identifying~$H$ with its
dual we obtain the continuous and dense embeddings
\[
H_1\subset H(\equiv H^*)\subset H^*_1.\vspace*{-2pt}
\]
It follows that
\[
_{H_1}\langle z,v\rangle_{H_1^*}=\langle z,v\rangle_H\qquad  \forall z\in
H_1,v\in H,\vspace*{-2pt}
\]
and that
$(H_1,H,H_1^*)$ is a~Gelfand triple.

The following is the main result of this section.\vspace*{-3pt}

\begin{thm}\label{thm6.2}
If $\alpha>0$, then $I_{K_\alpha}\in \operatorname{BV}(H,H_1)\cap
\mathbf{H}$.\vspace*{-3pt}
\end{thm}

\begin{pf}
First, for $\sigma_\alpha$ as in~(\ref{eq6.2}) we show that for each
$B\in\mathcal{B}(H)$ the function $r\mapsto\sigma_\alpha(r,B)$ is
in $H_1^*$ and that the map $B\mapsto\sigma_\alpha(\cdot,B)$ is in
fact an $H_1^*$-valued measure of bounded variation, that is,
\[
\sup\Biggl\{\sum_{n=1}^\infty\|\sigma_\alpha(\cdot,B_n)\|_{H_1^*}\dvtx
B_n\in\mathcal{B}(H),n\in\mathbb{N},H=\mathop{\dot{\bigcup}}_{n=1}^\infty
B_n\Biggr\}<\infty,\vspace*{-2pt}
\]
that is,
\begin{eqnarray*}
&&\sup\Biggl\{\sum_{n=1}^\infty\Biggl(\sum_{j=1}^\infty c_j^{-2}\biggl(\int_0^1\sigma
_\alpha(r,B_n)\sin(j\pi r)\,dr\biggr)^2\Biggr)^{1/2}\dvtx \\[-2pt]
&&
\hspace*{87pt} B_n\in\mathcal{B}(H),n\in
\mathbb{N},H=\mathop{\dot{\bigcup}}_{n=1}^\infty B_n\Biggr\}<\infty,\vspace*{-2pt}
\end{eqnarray*}
where $\dot{\bigcup}_{n=1}^\infty B_n$ means disjoint union.

For $\alpha>0$, we have
\begin{eqnarray*}
&&\sum_{n=1}^\infty\Biggl(\sum_{j=1}^\infty c_j^{-2}\biggl(\int
_0^1\sigma_\alpha(r,B_n)\sin(j\pi r)\,dr\biggr)^2\Biggr)^{1/2}\\[-2pt]
&&\qquad \leq\sum_{n=1}^\infty\Biggl(\sum_{j=1}^\infty c_j^{-2}\biggl(\int_0^1\sigma_\alpha
(r,B_n) \,dr\biggr)^2\Biggr)^{1/2}\\[-2pt]
&&\qquad \leq C\sum_{n=1}^\infty\int_0^1\sigma
_\alpha(r,B_n)\,dr\\[-2pt]
&&\qquad = C\int_0^1\sigma_\alpha(r)\,dr<\infty.
\end{eqnarray*}

Thus, $\sigma_\alpha$ in~(\ref{eq6.2}) is of bounded variation as an
$H_1^*$-valued measure. Hence by the theory of vector-valued measures
(cf.~\cite{2}, Section 2.1), there is a~unit vector field $n_\alpha
\dvtx H\rightarrow H_1^*$, such that $\sigma_\alpha=n_\alpha\|\sigma
_\alpha\|$, where
\[
\|\sigma_\alpha\|(B):=\sup\Biggl\{\sum_{n=1}^\infty\|
\sigma_\alpha(\cdot,B_n)\|_{H_1^*}\dvtx  B_n\in\mathcal{B}(H), n\in
\mathbb{N}, B=\mathop{\dot{\bigcup}}_{n=1}^\infty B_n\Biggr\}
\]
is a~nonnegative
measure, which is finite by the above proof. So~(\ref{eq6.2}) becomes
\[
\int_{K_\alpha}D^*(\varphi(\cdot)l)\,d\mu=\int{ }_{H_1} \langle
\varphi(x)l,n_\alpha(x)\rangle_{H_1^*}\|\sigma_\alpha\|(dx)
\qquad \forall l\in D(A), \varphi\in C_b^1(H),
\]
which by linearity extends to all $G\in(C_b^1)_{D(A)\cap H_1}$.
Thus by Theorem~\ref{thm3.1}(iii), we get that $I_{K_\alpha}\in \operatorname{BV}(H,H_1)$.

$I_{K_\alpha}\in\mathbf{H}$ follows by Remark~\ref{rem4.1}.
\end{pf}

\begin{remark}\label{rem6.3}
It has been proved by Guan Qingyang that $I_{K_\alpha
}$ is not in $\operatorname{BV}(H,H)$.
\end{remark}

\begin{thm}\label{thm6.4}
For $\alpha=0$, then there exist a~positive finite
measure $\|\sigma_0\|$ on~$H$ and a~Borel-measurable map
$n_0\dvtx H\rightarrow H_1^*$ such that $\|n_0(z)\|_{H_1^*}=1$ $\|
\sigma_0\|$-a.e., and
%
\begin{eqnarray}\label{eq6.4}
&&-\int_{K_0}\langle l,D\varphi\rangle\, d\nu-\int_{K_0}\varphi
(x)\langle x,l''\rangle\nu(dx)\nonumber\\[-8pt]\\[-8pt]
&&\qquad =\int{ }_{H_1} \langle\varphi
(x)l,n_0(x)\rangle_{H_1^*}\|\sigma_0\|(dx)\qquad  \forall l\in
D(A),\varphi\in C_b^1(H).\nonumber
\end{eqnarray}
\end{thm}

\begin{pf}
For $\alpha=0$ using that $|\sin(j\pi r)|\leq2j\pi r(1-r)
\ \forall r\in[0,1]$, we have
\begin{eqnarray*}
&&\sum_{n=1}^\infty\Biggl(\sum_{j=1}^\infty c_j^{-2}\biggl(\int
_0^1\sigma_0(r,B_n)\sin(j\pi r)\,dr\biggr)^2\Biggr)^{1/2}\\
&&\qquad \leq \sum_{n=1}^\infty\Biggl(\sum_{j=1}^\infty c_j^{-2}\biggl(\int_0^1\sigma_0(r,B_n)2j\pi
r(1-r)\,dr\biggr)^2\Biggr)^{1/2}\\
&&\qquad \leq C\sum_{n=1}^\infty\int_0^1\sigma
_0(r,B_n)r(1-r)\,dr\\[-2pt]
&&\qquad = C\int_0^1\sigma_0(r)r(1-r)\,dr<\infty.
\end{eqnarray*}
Thus, $\sigma_0$ in~(\ref{eq6.1}) is of bounded variation as an $H_1^*$-valued
measure. Hence by the theory of vector-valued measures (cf.~\cite{2}, Section
2.1), there is a~unit vector field $n_0\dvtx H\rightarrow H_1^*$, such that
$\sigma_0=n_0 \|\sigma_\alpha\|$, where
\[
\|\sigma_0\|(B):=\sup\Biggl\{
\sum_{n=1}^\infty\|\sigma_0(\cdot,B_n)\|_{H_1^*}\dvtx B_n\in\mathcal
{B}(H), n\in\mathbb{N}, B=\mathop{\dot{\bigcup}}_{n=1}^\infty B_n\Biggr\}
\]
is a~nonnegative measure, which is finite by the above proof. So the result
follows by~(\ref{eq6.1}).
\end{pf}

Since here $\mu(K_0)=0$, we have to change the reference measure of
the Dirichlet form. Consider
\[
\mathcal{E}^{K_0}(u,v)=\frac{1}{2}\int_{K_0}\langle Du,Dv\rangle
\,d\nu,\qquad u,v\in C_b^1(K_0).
\]
Since $I_{K_0}\in\mathbf{H}$ by Remark~\ref{rem4.1}, the closure of $(\mathcal
{E}^{I_{K_0}},C_b^1(K_0))$ is also a~quasi-regular local Dirichlet form
on $L^2(F;\rho\cdot\nu)$ in the sense of~\cite{17}, Chapter
IV, Definition 3.1. As
before, there exists a~diffusion process $M^{I_{K_0}}=(\Omega,\mathcal
{M},\allowbreak\{\mathcal{M}_t\},\theta_t,X_t,$ $P_z)$ on $F$ associated with
this Dirichlet form. $M^{I_{K_0}}$ will also be called distorted
OU-process on $K_0$. As before, $M^{I_{K_0}}$ is recurrent and
conservative. As before, we also have the associated PCAF and the Revuz
correspondence.

Combining these two cases: for $\alpha>0$ by Theorem~\ref{thm3.2} and for
$\alpha=0$ by the same argument as Theorem~\ref{thm3.2}, since we have~(\ref{eq6.4}),
we have the following theorem.

\begin{thm}\label{thm6.5}
Let $\rho:=I_{K_\alpha}, \alpha\geq0$ and
consider the measure $\|\sigma_\alpha\|$ and~$n_\alpha$ appearing in
Theorems~\ref{thm6.2} and~\ref{thm6.4}. Then there is an $\mathcal{E}^\rho
$-exceptional set $S\subset F$ such that $\forall z\in F\setminus S$,
under $P_z$ there exists an $\mathcal{M}_t$-cylindrical Wiener
process $W^z$, such that the sample paths of the associated distorted
OU-process~$M^\rho$ on $F$ satisfy the following: for $l\in D(A)$
%
\begin{eqnarray}\label{eq6.5}
\langle l,X_t-X_0\rangle&=&\int_0^t\langle l,dW_s\rangle+\frac
{1}{2}\int_0^t { }_{H_1} \langle l,n_\alpha(X_s)\rangle
_{H_1^*}\,dL_s^{\|\sigma_\alpha\|}\nonumber\\[-9pt]\\[-9pt]
&&{}-\int_0^t\langle Al, X_s\rangle\,ds\qquad  P_z\mbox{-a.e.}\nonumber
\end{eqnarray}
Here $L_t^{\|\sigma_\alpha\|}$ is the real valued PCAF associated
with $\|\sigma_\alpha\|$ by the Revuz correspondence with respect to
$M^\rho$, satisfying
%
\begin{equation}\label{eq6.6}
I_{\{X_s+\alpha\neq0\}}\,dL_s^{\|\sigma_\alpha\|}=0,\vadjust{\goodbreak}
\end{equation}
and for $l\in H_1$ with $l(r)\geq0$ we have
%
\begin{equation}\label{eq6.7}
\int_0^t { } _{H_1} \langle l,n_\alpha(X_s)\rangle_{H_1^*}\,dL_s^{\|
\sigma_\alpha\|}\geq0.
\end{equation}

Furthermore, for all $ z\in F$
%
\begin{equation}\label{eq6.8}
P_z\bigl[X_t\in C_0[0,1] \mbox{ for a.e. } t\in[0,\infty
)\bigr]=1.
\end{equation}
\end{thm}

\begin{pf}
For $\alpha>0$, the first part of the assertion follows by
Theorem~\ref{thm3.2} and the uniqueness part of Theorem~\ref{thm3.1}(ii). For $\alpha
=0$, the assertion follows by the same argument as in Theorem~\ref{thm3.2}.
(\ref{eq6.6}) and~(\ref{eq6.7}) follow by the property of $\sigma_\alpha$ in \cite
{28}. By~\cite{22}, p. 135, Theorem 2.4, we have $C_0[0,1]$ is a~Borel
subset of $L^2[0,1]$.
By~\cite{10}, (5.1.13), we have
\[
E_{\rho\mu}\biggl[\int_{k-1}^k1_{F\setminus C_0[0,1]}(X_s)\,ds\biggr]=\rho
\mu(F\setminus C_0[0,1])=0 \qquad  \forall k\in\mathbb{N},
\]
hence
\[
E_{\rho\mu}\biggl[\int_0^\infty1_{F\setminus C_0[0,1]}(X_s)\,ds\biggr]=0.
\]
Since $E_x[\int_0^\infty1_{F\setminus C_0[0,1]}(X_s)\,ds]$ is a~$0$-excessive function in $x\in K_\alpha$, it is finely continuous
with respect to the process $X$. Then for $\mathcal{E}^\rho$-q.e. $ z\in F$,
\[
E_z\biggl[\int_0^\infty1_{F\setminus C_0[0,1]}(X_s)\,ds\biggr]=0,
\]
thus, for $\mathcal{E}^\rho$-q.e. $z\in F$,
\[
P_z\biggl[\int_0^\infty1_{F\setminus C_0[0,1]}(X_s)\,ds=0\biggr]=1.
\]
As a~consequence, we have that $\Lambda_0:=\{X_t\in C_0[0,1]
\mbox{ for a.e. } t\in[0,\infty)\}$ is measurable and for
$\mathcal{E}^\rho$-q.e. $ z\in F$
\[
P_z(\Lambda_0)=1.
\]
As $\Lambda_0=\bigcap_{t\in\mathbb{Q}, t>0}\theta_t^{-1}\Lambda_0$
and since by~\cite{4} we have that the semigroup associated with $X_t$
is strong Feller, by the Markov property as in~\cite{8}, Lemma 7.1, we
obtain that for any $z\in F, t\in\mathbb{Q},t>0$,
\[
P_z(\theta_t^{-1}\Lambda_0)=1.
\]
Hence, for any $z\in F$ we have
\[
P_z\bigl[X_t\in C_0[0,1] \mbox{ for a.e. } t\in[0,\infty)\bigr]=1.
\]
\upqed
\end{pf}

\begin{remark}\label{rem6.6}
We emphasize that in the present situation it was
proved in~\cite{19}, Theorem 1.3, that for all initial conditions $x\in H$,
there exists a~unique strong solution to~(\ref{eq1.1}). By~\cite{28}, the
solution in~\cite{19} is associated to our Dirichlet form, hence\vadjust{\goodbreak}
satisfies~(\ref{eq6.5}) by Theorem~\ref{thm6.5}. Hence, it follows that the solution in
\cite{19}, Theorem 1.3, is solution to an infinite-dimensional Skorohod
problem.
\end{remark}

\section*{Acknowledgments}
The authors would like to thank Professor
Zhiming Ma and Qingyang Guan for their helpful discussions and
suggestions and also Giuseppe Da Prato and Ben Goldys for reaching an
earlier version of this paper.


%

\printaddresses

\end{document}